\documentclass[a4paper,pdftex,reqno,11pt]{amsart}

\usepackage{geometry}
\usepackage{graphicx}
\usepackage{courier}
\usepackage{times}
\usepackage{amsmath,amssymb}
\usepackage{hyperref}
\usepackage{color}

\newtheorem{Theorem}{Theorem}[section]
\newtheorem{Proposition}[Theorem]{Proposition}
\newtheorem{Conjecture}[Theorem]{Conjecture}
\newtheorem{Corollary}[Theorem]{Corollary}
\newtheorem{Lemma}[Theorem]{Lemma}

\newenvironment{Proof}
{\begin{trivlist}\item[]{{\sc Proof.}}}{\hfill{$\square$}\noindent\end{trivlist}}

\newtheorem{Definition}[Theorem]{Definition}

\newtheorem{Remark}[Theorem]{Remark}

\newcommand{\cv}[1]{\mathbf{\widetilde{{#1}}}}
\newcommand{\cardv}{\mathbf{\overline{{n}}}}
\newcommand{\irep}{integer representation}
\newcommand{\ireppt}{integer representation preserving types}
\newcommand{\ireppts}{integer representations preserving types}
\newcommand{\mirep}{minimum integer representation}
\newcommand{\mireppt}{minimum integer representation preserving types}
\newcommand{\msirep}{minimum sum integer representation}
\newcommand{\msireppt}{minimum sum integer representation preserving types}
\newcommand{\msireppts}{minimum sum integer representations preserving types}
\newcommand{\mwcs}{\mathcal{W}^m}
\newcommand{\mlcs}{\mathcal{L}^M}
\newcommand{\amount}{A}
\newcommand{\mwvs}{\mathcal{W}^{sm}}
\newcommand{\mlvs}{\mathcal{L}^{sM}}
\newcommand{\smin}{shift-minimal}
\newcommand{\smax}{shift-maximal}

\newcommand{\highlight}[1]{#1}

\begin{document}

\title{On minimum integer representations of weighted games}

\author{Josep Freixas$^a$}
\address{\textnormal{$^A$Department of Applied Mathematics III and Engineering School of Manresa, Technical University of Catalonia, Spain.}}
\author{Sascha Kurz$^b$}
\address{\textnormal{$^B$Fakult\"at f\"ur Mathematik, Physik und Informatik, Universit\"at Bayreuth, Germany,
email: sascha.kurz@uni-bayreuth.de, tel: +49 921 557352.}
}
\maketitle

{\small \noindent \textbf{Abstract:} We study {\mirep s} of weighted games, i.e.\ representations where the weights are integers and every other integer representation
  is at least as large in each component. Those {\mirep s}, if the exist at all, are linked with some solution concepts in game theory. Closing
  existing gaps in the literature, we prove that each weighted game with two types of voters admits a (unique)
  {\mirep}, and give new examples for more than two types of voters without a {\mirep}. We characterize the possible weights in {\mirep s} and give examples
  for $t\ge 4$ types of voters without a {\mireppt}, i.e.\ where we additionally require that the weights are equal within equivalence
  classes. }

{\small \medskip }



\section{Introduction}

\noindent
Simple games, or positive switching functions, can be viewed as models of voting systems in which a single alternative, such as a
bill or an amendment, is pitted against the status quo. Weighted games, or positive threshold functions, are possibly the most
interesting subclass of simple games. Roughly speaking, in a weighted game a non-negative weight $w_i$ is assigned to
each voter $1\le i\le n$ and a quota~$q$ is specified. As an abbreviation for a weighted game we use the
notation $[q\highlight{;}w_1,\dots,w_n]$. Winning coalitions are those that can force a victory, i.e.\ the sum of their weights equals or 
surpasses the quota.  Weighted games naturally appear in several different contexts apart from voting, like
reliability analysis of technical systems (see Ramamurthy \cite{Ram90}) or neural networks (see, among others, Elgot \cite{Elg60}
or Freixas and Molinero \cite{FrMo08}). 

The number of simple games on a fixed set $N=\{1,\dots,n\}$ is finite, of course, but it grows very rapidly with an increasing number
of voters~$n$ since we are dealing with sets of sets. Indeed, every family of pairwise independent subsets of $N$ can serve as the set
of minimal winning coalitions defining a simple game. Two subsets are independent  if neither contains the other. Families of independent
subsets are sometimes called ``Sperner families``, ``coherent systems``, or ``clutters``, and their enumeration and classification have
occupied mathematicians since Dedekind in the 19th century. In his 1897 work he determined the exact number of simple games with four or
fewer players. Since that time simple games have been investigated in a variety of different mathematical contexts. An account of some
of these works can be found in:  Sperner \cite{Spe28},  Isbell \cite{Isb59},  Golomb \cite{Gol59}, Muroga et
al. \cite{MTB70,MTK62}, Shapley and Shubik \cite{ShSh73}, Dubey and Shapley \cite{0409.90008}, Kurz and Tautenhahn \cite{min_weights,dedekind},
Freixas and Molinero \cite{FrMo09, FrMo10}, Krohn and Sudh\"olter \cite{0841.90134}, Keijzer et~al.\ \cite{keijzer,keijzer2}.
Although the number of weighted games compared with \highlight{the number of} simple games is small, it grows very rapidly and there 
\highlight{do not} exist enumeration results for more than nine voters.

Integer representations, i.e.\ where the weights $w_i$ and the quota $q$ are non-negative integers, are very common in practice
and {\mirep s}, if they exist, constitute the most efficient way to represent weighted games. Several algorithms to compute certain
power indices require \highlight{integer} weights and benefit from weights of small magnitude. When considering e.g.\ shareholders of a firm, 
\highlight{integer} weights,
i.e.\ the number of (equal) shares, arise naturally. Geometrically, the set of equivalent integer representations of a weighted game 
is an unbounded cone with or without a vertex. Hence, a natural question arises: For which weighted games
does a {\mirep} exist? Or, in other words, for which weighted games does the associated integer cone have a
vertex? Symmetric games, i.e.\ games where all players have an equivalent role in the game and, therefore, are
characterized by one single type of equivalent voters, admit a {\mirep}. But it is known that there does not always exist
a {\mirep} for a weighted game. Muroga et al.\ \cite{MTB70} in their exhaustive enumeration of threshold functions (or, equivalently,
weighted games) uncovered several cases with as few as eight players in which two symmetric players must be given different weights in
a {\msirep}; e.g.\ $[12\highlight{;}7,6,6,4,4,4,3,2]=[12\highlight{;}7,6,6,4,4,4,2,3]$. Here a {\msirep} is an {\irep} such that the
sum of weights $\sum_{i=1}^n w_i$ is minimal. Moreover, they verified that all weighted games with less than eight players admit a
{\mirep}. We can easily check that this example consists of four types of players (a type here is an element of a
partition of $N$ formed by equivalent voters), and each type contains players with the same weights except for the last
type, which contains players with weights $3$ and $2$.

To our knowledge it is not known whether there exist weighted games without a
{\mirep} with either two or three types of players. The main goal of this paper is to ascertain what occurs for these two
cases, filling the existing gap in the theory of weighted games. Previous to Muroga et~al's example, Isbell \cite{Isb59}
had exhibited a remarkable $12$-player example in which the affected players are not symmetric. Thus, even if we additionally require
that all players of equal type have equal weights, the existence of a {\mireppt} is not guaranteed. 
Freixas and Molinero \cite{FrMo09, FrMo10} uncovered several cases of weighted games without a {\mireppt} with as few as
$9$~players and checked the nonexistence of such examples for less than $9$ players; see also \cite{min_weights}. All the examples
they listed have at least $5$~types of players. So quite naturally, we want to ascertain what occurs for less than
$5$~types. We would like to remark that homogeneous games\footnote{A weighted game is called \textit{homogeneous} if it admits
a representation where all minimal winning coalitions have the same weight.} admit a {\mirep} as shown by Ostmann~\cite{0628.90100}.

A natural third issue emerges to be significant, whenever there does not exist a {\mirep} for a weighted game
(either preserving types or not). \highlight{In that situation} at least two integer representations are minimal, but is it possible to generate
weighted games with more than two \highlight{minimal  representations?} \highlight{Since integer representations which attain the minimum
possible sum of weights are minimal, we ask more generally for constructions of weighted games with an arbitrary
number of {\msirep s}.} As far as we know, all the previously published examples without a {\mirep} (either preserving types or not)  have
only two {\msirep s}. Additional results\highlight{,} we introduce here\highlight{,} comprise: bounds on the number of non-isomorphic weighted games as
a function of the number of voters and the number of types of voters, and the existence of a weighted game in {\mirep} for any pair of two
coprime integer weights. 

Minimum integer representations of weighted games are important in game theory: Peleg~\cite{peleg} proved that for \textit{homogeneous} weighted
\textit{decisive} games the nucleolus (a well-known solution concept in game theory) coincides with the {\mireppt}. Also, in the cases where there
is no {\mireppt}, there are connections linking a {\msireppt} with the least core (another solution concept) and
the nucleolus of weighted decisive games \cite{0841.90134}. 

The remainder of the paper is organized as follows. In Section~\ref{sec_basics} we precisely define the classes of complete simple games and weighted games.
For complete simple games we state a parameterization theorem by Carreras and Freixas in Subsection~\ref{subsec_parameterization_csg}, which
completely characterizes these objects up to isomorphism using linear inequalities. The subclass of weighted games can be defined via the
non-emptiness of certain polytopes as outlined in Subsection~\ref{subsec_weighted}. The details on {\mirep s} are stated in Subsection~\ref{subsec_minimum_integer_representation}. In Section~\ref{sec_3} we present constructions for weighted games without
a {\mirep} for small $t$ (Subsection~\ref{subsec_3_1}) and for those with more than two {\msirep s}
(Subsection~\ref{subsec_3_2}). In Subsection~\ref{subsec_3_3} we study the question of which weights may occur in a {\mirep}. Our
main theorem, that each weighted game with two types of voters admits a {\mirep}, is given in Section~\ref{sec_4}. Implications
for the enumeration or bounds on the number of weighted games, which arise as a byproduct of our previous results, are
briefly stated in Section~\ref{sec_5}. We end with a conclusion in Section~\ref{sec_conclusion}.

\section{Simple games, complete simple games and weighted games}
\label{sec_basics}

\noindent
From a more general point of view, binary voting systems, i.e.\ those where each voter has the option to vote \textit{yes}
or \textit{no}, which then is condensed by a certain voting rule, can be represented by a characteristic function $\chi:2^N\rightarrow\{0,1\}$,
where $N:=\{1,\dots,n\}$ is the set of voters and $2^N$ denotes the set $\{U\mid U\subseteq N\}$ of all subsets of $N$. A quite natural
monotonicity assumption on $\chi$ leads to a very prominent class of binary voting systems.

\begin{Definition}
A \textbf{simple game} is a function $\chi:2^N\rightarrow\{0,1\}$, which satisfies $\chi(\emptyset)=0$,
$\chi(N)=1$, and $\chi(U')\le\chi(U)$ for all $U'\subseteq U\subseteq N$, where $N$ is a finite set.
\end{Definition}

So, if we identify $2^N$ with $\{0,1\}^n$, each simple game is a monotone Boolean function and except for the all-zero function and the all-one function
all monotone Boolean functions are simple games. We will call a subset $U\subseteq N$ a coalition.

\begin{Definition}
A coalition $U\subseteq N$ of a simple game $\chi$ is called \textbf{winning} if $\chi(U)=1$ and \textbf{losing} otherwise. 
A coalition $U$ is called a \textbf{minimal winning} coalition if $\chi(U)=1$ and $\chi(U')=0$ for all
proper subsets $U'$ of $U$. Similarly, a coalition $U$ is called a \textbf{maximal losing} coalition if $\chi(U)=0$ and $\chi(U')=1$ for
all proper supersets $U'$ of $U$. By $\mathcal{W}$ we denote the set of winning coalitions and by $\mathcal{L}$ the set of losing
coalitions for a given simple game. The restrictions to minimal winning coalitions and maximal losing coalitions are denoted by
$\mwcs$ and $\mlcs$, respectively.
\end{Definition}

We have $\mathcal{W}\cup\mathcal{L}=2^N$ and remark that either $\mwcs$ or $\mlcs$
uniquely characterizes a simple game; see e.g.\ \cite{0943.91005} for the details and additional facts on simple games. 
A well studied subclass of simple games (and superclass of weighted games) arises from Isbell's desirability relation
\cite{0083.14301}:
\begin{Definition}
We write $i\sqsupset j$ (or $j \sqsubset i$) for two voters $i,j\in N$ if we have $\chi\Big(\{i\}\cup U\backslash\{j\}\Big)\ge\chi (U)$
for all $\{j\}\subseteq U\subseteq N\backslash\{i\}$ and we abbreviate $i\sqsupset j$, $j\sqsupset i$ by $i~\square~j$. A simple game $\chi$
is called \textbf{complete simple game} (also called a ``directed game", see \cite{0841.90134}, \highlight{or a ``linear game", see \cite{0943.91005}})
if the binary relation $\sqsupset$ is a total preorder, i.e.\
\begin{itemize}
  \item[(1)] $i\sqsupset i$ for all $i\in N$,
  \item[(2)] $i\sqsupset j$ or $j\sqsupset i$ for all $i,j\in N$, and
  \item[(3)] $i\sqsupset j$, $j\sqsupset h$ implies $i\sqsupset h$ for all $i,j,h\in N$.
\end{itemize}
\end{Definition}

W.l.o.g.\ we assume $1\sqsupset 2 \sqsupset \dots \sqsupset n$ in the following. Whenever $i~\square~j$, voter~$i$
is as influential in the game as voter~$j$, meaning that it does not matter which one of both takes part in a coalition, i.e.\ the status of
the coalition (winning or losing) does not change after a swap of two equally desirable voters. We can partition
the whole set~$N$ of voters into equivalence classes $N_1,\dots, N_t$ and say that the complete simple game consists of $t$~types of voters. 
By $n_i$ we denote the cardinality of the set $N_i$ for $1\le i\le t$. Coalitions are categorized into different types,
which can be described by a vector $(m_1,\dots,m_t)$ meaning $m_i$-out-of-$n_i$ voters (from the set $N_i$) for $1\le i\le t$. 

Let us consider an example with $n_1=n_2=2$. Due to the assumed ordering of the players we have $N_1=\{1,2\}$ and $N_2=\{3,4\}$. With this
the vector $(1,1)$ is the type of the coalitions $\{1,3\}$, $\{1,4\}$, $\{2,3\}$, and $\{2,4\}$. Since we have $1~\square~2$ and $3~\square~4$ either all those
four coalitions are winning or they are all losing and we can therefore speak of a winning or a losing vector.

\begin{Definition}
\label{def_winning_vector}
Let $\chi$ be a simple game and $N_h$, \highlight{$1\le h\le t$,} be the classes of equally desirable voters. We call a
vector $\cv{m}:=(m_1,\dots,m_t)$, where $0\le m_h\le \left|N_h\right|$ for $1\le h\le t$, a
\textbf{winning vector} if $\chi(U)=1$, where $U$ is an arbitrary subset of $N$ containing exactly $m_h$ elements of $N_h$
for $1\le h\le t$. Analogously, we call such a vector a \textbf{losing vector} if $\chi(U)=0$, where $U$ is an arbitrary
subset of $N$ containing exactly $m_h$ elements of $N_h$ for $1\le h\le t$.
\end{Definition}

In the following we will always use a tilde and bold notation to indicate a vector representing a type
of a coalition. The concept of inclusion has to be slightly modified for vectors, i.e.\ types of coalitions:

\begin{Definition}
  \label{def_smaller_vector}
  For two vectors $\cv{a}=(a_1,\dots,a_t)$ and $\cv{b}=(b_1,\dots,b_t)$, representing types of coalitions in a complete
  simple game, we write $\cv{a}\preceq \cv{b}$ if we have
  $
    \sum\limits_{i=1}^{k} a_i \le \sum\limits_{i=1}^{k} b_i
  $
  for all $1\le k\le t$. For $\cv{a}\preceq \cv{b}$ and $\cv{a}\neq \cv{b}$ we use
  $\cv{a}\prec\cv{b}$ as an abbreviation and say that they are comparable vectors with vector $\cv{a}$ being smaller
  than vector $\cv{b}$. If neither $\cv{a}\preceq \cv{b}$ nor $\cv{b}\preceq \cv{a}$ holds, we write $\cv{a}\bowtie \cv{b}$ 
  and say that vector $\cv{a}$ and vector $\cv{b}$ are incomparable.
\end{Definition}
If $(1,1)$ is a winning vector in our example, so is $(2,0)$ while nothing can be deduced for vector $(0,2)$. 
With Definition \ref{def_smaller_vector} at hand, we can define:
\begin{Definition}
\label{def_shift_minimal}
A vector $\cv{m}=(m_1,\dots,m_t)$ in a complete simple game $\chi$ is a \textbf{\highlight{\smin} winning vector}
if $\cv{m}$ is a winning vector and every vector \highlight{$\cv{m}\mathbf{'}$} with $\cv{m}\mathbf{'}\prec\cv{m}$ is losing. Analogously,
a vector $\cv{m}$ is a \textbf{\highlight{\smax} losing vector} if $\cv{m}$ is a losing vector and every vector
\highlight{$\cv{m}\mathbf{'}$ with} $\cv{m}\mathbf{'}\succ\cv{m}$ is winning.
\end{Definition}

Similarly as for simple games, where the set $\mwcs$ or $\mlcs$ with the inclusion are enough to generate the entire
set of winning coalitions $\mathcal{W}$, for complete simple games the sets  
$\mwvs$ and $\mlvs$ of the \highlight{{\smin}} winning vectors (representing types of coalitions) and the maximal losing vectors uniquely
characterize the complete simple
game with the operation $\succeq$. Weighted games, which are a subclass of complete simple games, are now formally introduced as follows:
\begin{Definition}
  \label{def_weighted_voting_game}
  A simple game $\chi$ is called a \textbf{weighted game} (or simply \textbf{weighted}) if there exists a quota
  $q\in\mathbb{R}_{>0}$ and weights $w_1,\dots,w_n\in\mathbb{R}_{\ge 0}$ such that $\chi(U)=1$ \emph{if and only if} $\sum_{i\in U} w_i\ge q$.
  As an abbreviation we utilize the notation $\chi=[q\highlight{;}w_1,\dots,w_n]$ or simply $\chi=[q\highlight{;}w]$ whenever the weight
  vector $w=(w_1,\dots,w_n)$ is specified.
\end{Definition}

As an example we consider the weighted game $[4\highlight{;}3,2,1,1]$ (which is the same as $[3\highlight{;}2,1,1,1]$), where we have
$1\sqsupset 2~\square~3~\square~4$ for the voters, i.e.\ $n_1=1$ and $n_2=3$. The \highlight{\smin} winning 
vectors are given by $(1,1)$, $(0,3)$ and the \highlight{smax} losing vectors
are given by $(1,0)$, $(0,2)$. Since $(1,2)\succ(1,1)$ the coalition type $(1,2)$ is also winning and $(0,2)$ is
losing due to $(0,2)\prec(0,3)$. For a more extensive overview on binary voting methods
we refer the interested reader to~\cite{0943.91005}.

\subsection{A parameterization theorem for complete simple games}
\label{subsec_parameterization_csg}
Carreras and Freixas have given a full parameterization of complete simple games
in  \cite{complete_simple_games}. To this end we denote the (decreasing) lexicographic \highlight{(strict)} order by $\gtrdot$, i.e.\
we have  $(a_1,\dots,a_n) \gtrdot (b_1,\dots,b_n)$ iff there is an index $1\le h\le n$ with $a_i=b_i$ for all $1\le i<h$ and $a_h>b_h$.
An example is given by $(1,2,1)\gtrdot(1,1,3)$.

\begin{Theorem}(Carreras and Freixas, 1996)
  \label{thm_characterization_cs}

  \vspace*{0mm}

  \noindent
  \begin{itemize}
   \item[(a)] \highlight{Consider a vector} $$\cardv=(n_1,\dots,n_t)\in\mathbb{N}_{>0}^t$$ and \highlight{a} matrix
              $$\mathcal{M}=\begin{pmatrix}m_{1,1}&m_{1,2}&\dots&m_{1,t}\\m_{2,1}&m_{2,2}&\dots&m_{2,t}\\
              \vdots&\vdots&\ddots&\vdots\\m_{r,1}&m_{r,2}&\dots&m_{r,t}\end{pmatrix}=
              \begin{pmatrix}\cv{m}_1\\\cv{m}_2\\\vdots\\\cv{m}_r\end{pmatrix}.$$
              \highlight{If they satisfy} the following properties:
              \begin{itemize}
               \item[(i)]   $m_{1,1}>0$, $0\le m_{i,j}\le n_j$, $m_{i,j}\in\highlight{\mathbb{N}_{\ge 0}}$ for $1\le i\le r$, $1\le j\le t$,
               \item[(ii)]  $\cv{m}_i\bowtie\cv{m}_j$ for all $1\le i<j\le r$,
               \item[(iii)] for each $1\le j<t$ there is at least one row-index $i$ such that
                            $m_{i,j}>0$, $m_{i,j+1}<n_{j+1}$, and
               \item[(iv)]  $\cv{m}_i\gtrdot \cv{m}_{i+1}$ for $1\le i<t$,
              \end{itemize}
              \highlight{then} there exists a complete simple game $\chi$ associated to $\left(\cardv,\mathcal{M}\right)$ with $\cardv$ as 
              a vector of the cardinalities of the equivalence classes and matrix $\mathcal{M}$, where the rows consist of the \highlight{\smin}
              winning vectors.
   \item[(b)] Two complete games $\left(\cardv_1,\mathcal{M}_1\right)$ and $\left(\cardv_2,\mathcal{M}_2\right)$
              are isomorphic \emph{if and only if} $\cardv_1=\cardv_2$ and $\mathcal{M}_1=\mathcal{M}_2$.
  \end{itemize}
\end{Theorem}

\noindent
In such a vector/matrix representation of a complete simple game the number of voters $n$ is determined by $n=\sum_{i=1}^t n_i$.
Although Theorem~\ref{thm_characterization_cs} looks technical at first glance, the necessity of the required properties can be
explained easily. \highlight{First we observe that} $n_j\ge 1$, $m_{1,1}>0$, and $0\le m_{i,j}\le n_j$ must hold for $1\le i\le r$, $1\le j\le t$. If
$\cv{m}_i\preceq\cv{m}_j$ or $\cv{m}_i\succeq\cv{m}_j$ then we would have $\cv{m}_i=\cv{m}_j$ or either $\cv{m}_i$ or $\cv{m}_j$
cannot be a \highlight{\smin} winning vector. If for a column-index $1\le j<t$ we have $m_{i,j}=0$ or
$m_{i,j+1}=n_{j+1}$ for all $1\le i\le r$, then we can check that $g~\square~h$ for all $g\in N_j$, $h\in N_{j+1}$,
which is a contradiction to the definition of the classes $N_j$ and therefore also for the numbers $n_j$. \highlight{A} complete simple
game does not change if two rows of the matrix $\mathcal{M}$ are interchanged. Thus we must require some specific ordering
of the rows to avoid duplicities, e.g.\ $\gtrdot$.

\highlight{If all voters are equivalent, i.e.\ $t=1$, there is a unique \highlight{\smin} winning vector, i.e.\ $r=1$. In this case}
the requirements of Theorem~\ref{thm_characterization_cs} are reduced to
$1\le m_{1,1}\le n_1=n$. Also for $t=2$ one can easily give a more
compact formulation for the requirements in Theorem \ref{thm_characterization_cs}. A complete description of the possible values
$n_1,n_2,m_{1,1},m_{1,2}$ corresponding to a complete simple game with parameters $n$, $t=2$, and $r=1$ is given by
\begin{equation}
  1\le n_1\le n-1,\quad
  n_1+n_2=n,\quad
  1\le m_{1,1} \le n_1, \,\,\text{ and }\,\,
  0\le m_{1,2}\le n_2-1.\label{compact_ilp_2_1}
\end{equation}
For $t=2$ and $r\ge 2$ such a complete and compact description is given by
\begin{equation}
  1\le n_1\le n-1,\,\,
  n_1+n_2 = n,\text{ and }
  m_{i,1}\ge m_{i+1,1}+1, \,\,
  m_{i,1}+m_{i,2}+1\le m_{i+1,1}+m_{i+1,2} \label{compact_ilp_2_ge_2}
\end{equation}
for all $1\le i\le r-1$.

\subsection{Recognizing and representing weighted games}
\label{subsec_weighted}

\noindent
In Definition~\ref{def_weighted_voting_game} we have introduced the notation $[q\highlight{;}w_1,\dots,w_n]$,
consisting of a quota $q$ and weights $w_i$, for a weighted game. As mentioned in the introduction there are several
representations for the same weighted game,
e.g.\ $[3\highlight{;}2,1,1,1]$, $[4\highlight{;}3,2,1,1]$, $[11\highlight{;}9,5,5,4]$, $[q\highlight{;}q-1,x,x,x]$ and $[q\highlight{;}q-2,x,x,x]$
with $q\ge 6$ and $\left\lceil\frac{q}{3}\right\rceil\le x \le\left\lfloor\frac{q-1}{2}\right\rfloor$ all represent the same weighted
game because the subsets of $N$ whose weights equal or surpass the quota are invariant for all of them.

So in order to check whether two weighted games are equivalent, it makes sense to have a
closer look at the underlying discrete structure as a simple game, i.e.\ its characteristic function
$\chi:2^N=\left\{U\mid U\subseteq N\right\}\rightarrow\{0,1\}$. As weighted games are complete
simple games we often find it useful to represent the game using the matrix representation of the previous
subsection, especially if we use different weighted representations for the same game or different weights
within an equivalence class of voters.

To decide whether a given complete simple game is weighted, we can utilize a linear program; see \cite{0943.91005} for an
overview on other methods. From Definition~\ref{def_weighted_voting_game} and the notion of minimal
winning and maximal losing coalitions we can conclude that a simple game is weighted \emph{if and only if} the
following system of linear inequalities is feasible:
\begin{equation}
  \sum_{i\in S} w_i\ge q\,\, \forall S\in\mwcs,\,\,\sum_{i\in T} w_i< q\,\, \forall T\in\mlcs,\,\,
  q\in\mathbb{R}_{>0},\text{ and }w_i\in\mathbb{R}_{\ge 0}\,\,\forall 1\le i\le n.\label{ie_feas0}
\end{equation}
As strict inequalities, i.e., $<$ or $>$, might lead to ill-defined optimization problems like e.g.\ maximize $x$
subject to $x < 1$, we use an equivalent formulation instead:
\begin{equation}
  \sum_{i\in S} w_i\ge q\,\, \forall S\in\mwcs,\,\,\sum_{i\in T} w_i\le q-1\,\, \forall T\in\mlcs,
  \text{ and }w_i\ge 0\,\,\forall 1\le i\le n.\label{ie_feas1}
\end{equation}
As $\mlcs$ is not empty and the $w_i$ are non-negative, the inequality $q>0$ is implied by
$\sum_{i\in T} w_i\le q-1$. By rescaling the weights we may achieve that the difference $q-\max_{T\in\mlcs}\sum_{i\in T} w_i$
is as large as desired, e.g.\ at least $1$. Of course here we already have integer representations in mind, i.e.\ where
we additionally request $w_i\in\highlight{\mathbb{N}_{\ge 0}}$  (see Definition~\ref{def_integer_representation}). 
The fact that each weighted game is also a complete simple game can be used to reduce inequality
system~(\ref{ie_feas1}). 
\begin{Lemma}
  \label{lemma_ie_feas2}
  Given a complete simple game $\chi$ with $t$ equivalence classes of voters the inequality system~(\ref{ie_feas1}) has
  a solution \emph{if and only if}
  \begin{equation}
    \cv{x}^T w\ge q\,\, \forall\, \cv{x}\in\mwvs,\,\,\cv{y}^Tw\le q-1\,\, \forall\, \cv{y}\in\mlvs,
    w_i\ge w_{i+1}+1\,\,\forall 1\le i\le t-1, \text{ and } w_t\ge 0.\label{ie_feas2}
  \end{equation}
  has a solution.
\end{Lemma}
\begin{Proof}
  Let us at first assume that $(q,w)$ is a feasible solution of (\ref{ie_feas2}). By setting $q'=q$ and $w_i'=w_h\highlight{\ge 0}$ for all
  $i\in N_h$ we will obtain a feasible solution $(q',w')$ for (\ref{ie_feas1}). Now let $S\in\mwcs$ be a minimal winning coalition,
  $\cv{x}'$ its corresponding type, \highlight{and let} $\cv{x}\in\mwvs$ \highlight{be a} vector with $\cv{x}\preceq\cv{x}'$.
  With this we have
  $$
    \sum_{i\in S} w_i'=\cv{x}'^Tw\ge \cv{x}^Tw\ge q
  $$
  due to $w_1>w_2>\dots >w_t\ge 0$ for all $1\le i\le t$. Similarly, for a maximal losing coalition $R\in\mlcs$ with corresponding
  type $\cv{y}'$, let $\cv{y}\in\mlvs$ be a vector with $\cv{y}\succeq\cv{y}'$, so that
  $$
    \sum_{i\in R} w_i'=\cv{y}'^Tw\le \cv{y}^Tw\le q-1.
  $$
  
  \medskip
  
  For the other direction let $(q',w')$ be a feasible solution of~(\ref{ie_feas1}). 
  \highlight{One can easily} check
  that $(q',w'')$, where $w_i''=\frac{\sum_{j\in N_i} w_j'}{|N_i|}$ for all $1\le i\le t$, is a feasible solution of~(\ref{ie_feas2}).
  %
  %
\end{Proof}

We would like to remark that those complete simple games which are not weighted can be represented as a finite intersection
of weighted games, a construction which is  also used in practice \cite{pre05664394}.

\subsection{Minimum integer representations}
\label{subsec_minimum_integer_representation}

\noindent
In the previous section we have already seen some different representations of weighted games, e.g.\ we may assume that the difference
between the weight of a winning coalition and the weight of a losing coalition is at least one. A special kind of representation restricts the quota
and the weights to integers:

\begin{Definition}
  \label{def_integer_representation}
  For a given weighted game $\chi$, with minimal winning coalitions $\mwcs$, maximal losing coalitions $\mlcs$, and
  $t$ equivalence classes of voters, a vector $(q,w_1,\dots,w_n)\in\highlight{\mathbb{N}_{\ge 0}^{n+1}}$ is called an \textbf{\irep} if
  it is a feasible solution of Inequality system~(\ref{ie_feas1}). If we have $w_i=w_j$ for all $i,j\in N_h$, where $1\le h\le t$, then
  we speak of an \textbf{\ireppt}.
\end{Definition}

We remark that each feasible solution $(q,w_1,\dots,w_n)\in\highlight{\mathbb{N}_{\ge 0}^{n+1}}$ of (\ref{ie_feas0}) also satisfies
inequality system~(\ref{ie_feas1}). 
Given an {\irep} we can easily construct a (possibly non-integer)
representation, where the weights are equal within equivalence classes of voters by averaging the weights in each equivalence
class, as done in the proof of Lemma~\ref{lemma_ie_feas2}. (Every convex combination of solutions of an LP is itself a solution.)

\begin{Definition}
  Given an {\irep} $(q,w_1,\dots,w_n)$ for a weighted game $\chi$ with equivalence classes $N_1,\dots,N_t$ the 
  \textbf{averaged representation} $(q,w_1',\dots,w_t')$ is given by $w_h'=\frac{\sum_{i\in N_h} w_i}{|N_h|}$.
\end{Definition}
So indeed each weighted game admits an {\ireppt}.

\begin{Definition}
  Given a weighted game $\chi$ we call an {\irep} $(q,w_1,\dots,w_n)$ a \textbf{\msirep}, if
  we have $\sum_{i=1}^n w_i\le \sum_{i=1}^n w_i'$ for all {\irep s} $(q',w_1',\dots,w_n')$. Similarly we call an
  {\irep} $(q,w_1,\dots,w_n)$ preserving types a \textbf{\msireppt}, if
  we have $\sum_{i=1}^n w_i\le \sum_{i=1}^n w_i'$ for all {\irep s} $(q',w_1',\dots,w_n')$ preserving types.
\end{Definition}
We remark that each weighted game admits a {\msirep} and a {\msireppt}, but there can exist several
such representations. Introducing integer variables changes the linear programs (\ref{ie_feas1}) and
(\ref{ie_feas2}) to integer linear programs (ILP), whose solution is $NP$-hard in general.
So, if we minimize the sum of weights $\sum_{i=1}^n w_i$ subject to the constraints in inequality
system~(\ref{ie_feas1}) restricted to integer variables, each optimal solution corresponds to a {\msirep}. Similarly, if
we minimize the sum of weights $\sum_{i=1}^t n_iw_i$ subject to the constraints in inequality system~(\ref{ie_feas2}) restricted
to integer variables, each optimal solution corresponds to a {\msireppt}. To our knowledge there is no known polynomial time
algorithm to determine a {\msirep}. For some algebraic techniques, to determine a {\msirep}, we refer 
the interested reader to \cite{algebraic}. 

By considering the following LP-relaxation of the ILP for the value of a {\msirep} we can obtain a reasonable lower bound \highlight{for
the sum of weights in an {\msirep}}:
\begin{eqnarray}
  \!\!\!\!\!\!&&\min \sum_{i=1}^t w_in_i\label{lp_minsum}\\
  \!\!\!\!\!\!&&\text{s.t.}\  \cv{x}^T w\ge q\,\, \forall\, \cv{x}\in\mwvs,\,\,\cv{y}^Tw\le q-1\,\, \forall\, \cv{y}\in\mlvs,
    w_i\ge w_{i+1}+1\,\,\forall 1\le i\le t-1, \text{ and } w_t\ge 0.\nonumber
\end{eqnarray}

\begin{Lemma}
  \label{lemma_lb_msirep}
  For a given weighted game $\chi$ with $t$ equivalence classes of voters let $\varphi$ be the optimal target value of the minimization
  problem~(\ref{lp_minsum}), then we have $\sum_{i=1}^n w_i' \ge \varphi$ for all {\irep s} $(q',w_1',\dots,w_n')$ of $\chi$.
\end{Lemma}
\begin{Proof}
  \highlight{For a given {\irep} $(q',w_1',\dots,w_n')$ we show that the averaged representation $w_i=\frac{\sum_{j\in N_i} w_j'}{|N_i|}$,  $q=q'$ 
  is a feasible solution of inequality system~(\ref{ie_feas2}) attaining the same sum of its  weights as the initial {\irep}.}
  
  As in the proof of Lemma~\ref{lemma_ie_feas2} we have $w_{j_1}'>w_{j_2}'$ for all $j_1\in N_i$, $j_2\in N_{i+1}$. Since the $w_j'$
  are integers we conclude $w_{j_1}'\ge w_{j_2}'+1$ so that $w_i\ge w_{i+1}+1$ for all $1\le i\le t-1$. 
\end{Proof}

A more restrictive {\irep} asks for the minimum possible weight for each player simultaneously:
\begin{Definition}
  \label{def_minimum_integer_representation}
  An {\irep} $(q,w_1,\dots,w_n)$ for a weighted game $\chi$ is called \textbf{\mirep} if for
  all {\irep s} $(q',w_1',\dots,w_n')$ of $\chi$ we have $w_i\le w_i'$ for all $1\le i\le n$. If we restrict the allowed
  representations to those where the voters of the same equivalence class $N_i$ have an equal weight, we speak of a
  \textbf{\mireppt}. 
\end{Definition}
In other words, a {\mirep}, if it exists, is the least element in the partial order of component-wise
comparison of the feasible weight vectors.


In general, both representations need not exist and indeed in this paper we study conditions where they exist and give
examples where they do not exist. We would like to note that \highlight{each} {\mirep} for a weighted game \highlight{is also} a
{\mireppt}, \highlight{since otherwise the weights could be permuted within equivalence classes of voters.} \highlight{On the other
hand, the existence of a {\mireppt} does not imply the existence of a {\mirep}.} The example
$[12\highlight{;}7,6,6,4,4,4,3,2]=[12\highlight{;}7,6,6,4,4,4,2,3]$ from the introduction has $(14,8,7,7,5,5,5,3,3)$ as a {\mireppt}.

\section{Generating conspicuous examples of games without a {\mirep}}
\label{sec_3}

\noindent
Motivated by the existence of weighted games without a {\mirep} for more than
three equivalence classes of voters; see e.g.\ Table~3 and Table~4 of \cite{FrMo09}, we are concerned in this section with this problem in the special
case of $t=3$~types of voters. As we shall see below, we propose a procedure to generate weighted games with three types of
voters without a {\mirep} in Subsection~\ref{subsec_3_1} based on the famous Coin-Exchange Problem of Frobenius \cite{1114.52013}.
Similarly, the existence of weighted games without a {\mireppt} is known for more than four equivalence
classes of voters; see e.g.\ Table~2 in \cite{FrMo10}. Thus the case $t=4$ is
under study here and we also propose a procedure to generate weighted games with four types of voters without a
{\mireppt} in Subsection~\ref{subsec_3_1}. Another objective of this section is to generate examples of weighted games
with more than two {\msirep s}, which is outlined in Subsection~\ref{subsec_3_2}. Finally, Subsection~\ref{subsec_3_3}
concerns weighted games with a {\mirep} of coprime weights. 

The Coin-Exchange Problem of Frobenius considers
$n\ge 2$ integers $0<a_1<\dots<a_n$ with $\gcd(a_1,\dots,a_n)=1$ as denominations of $n$~different coins. We say
that a certain amount of money~$\amount\in\highlight{\mathbb{N}_{\ge 0}}$ can be \textbf{represented} by the given coins, if there are $n$ numbers
$x_i\in\highlight{\mathbb{N}_{\ge 0}}$ such that $\amount=\sum\limits_{i=1}^n a_ix_i$. As an abbreviation we denote the set of representable
integers $\amount$ by $\langle a_1,\dots,a_n\rangle$.

If $a_1>1$ then some $\amount$ cannot be represented, e.g.\ there do no exist representations for all $\amount\in\{1,\dots,a_1-1\}$.
The largest such $\amount$ for a given problem is called the Frobenius number $g(a_1,\dots,a_n)$. Well-known results in this context are
$g(a_1,a_2)=(a_1-1)(a_2-1)-1$ and that exactly $\frac{1}{2}(g(a_1,a_2)+1)=\frac{1}{2}(a_1-1)(a_2-1)$ non-negative
integers are not representable for $\gcd(a_1,a_2)=1$. As an example we consider $a_1=3$, $a_2=5$, where the set of
non-negative integers which are  not representable is given by $\highlight{\mathbb{N}_{\ge 0}} \backslash \langle a_1,a_2\rangle=\{1,2,4,7\}$.

Almost all of the following constructions contain the game
$\chi_{a,b}=[ab\highlight{;}\overset{a}{\overbrace{b,\dots,b}},\overset{b}{\overbrace{a,\dots,a}}]$, where $b>a\ge 1$ are coprime integers, as a subgame,
i.e.\ the winning coalitions of $\chi_{a,b}$ are winning coalitions in the larger game and similarly the losing coalitions of $\chi_{a,b}$
are losing coalitions of the larger game. Our first aim is to prove a lower bound on the sum of weights of a {\msirep} of $\chi_{a,b}$. To
this end we utilize B\'ezout's identity stating that there exist integers $u,v\in\mathbb{Z}$ with $ua+bv=\gcd(a,b)=1$, which can
be computed using the extended Euclidean algorithm.

\begin{Lemma}
  \label{lemma_xgcd}
  For coprime integers $b>a\ge 1$ there exist $u,v\in\mathbb{N}_{>0}$ with $ub-va=1$, $u\le a$, $v<b$.
\end{Lemma}

\begin{Lemma}
  \label{lemma_gcd_special}
  For coprime integers $b>a\ge 1$ there exist $u,v\in\mathbb{N}_{\ge 0}$ with $ub+va=ab-1$, $0\le u\le a-1$, and $1\le v\le b-1$. 
\end{Lemma}
\begin{Proof}
  Using Lemma~\ref{lemma_xgcd} and the identity $(a-u)b+va=ab-1$ yields the stated result.
\end{Proof}

\highlight{In the following we will often use the existence of those integers $u,v$ without explicitly referring to 
Lemma~\ref{lemma_gcd_special}.} We remark that the (unique) existence of such a pair $(u,v)$ of integers can be concluded from Popoviciu's theorem,
which counts the number of representations for a given amount $N$ using two coprime integer coins $a$ and $b$.

\begin{Lemma}
  \label{lemma_subgame_lb}
  For every {\irep} $(q,w_1,\dots,w_n)$ of $\chi_{a,b}$ we have $\sum_{i=1}^n w_i\ge 2ab$.
\end{Lemma}
\begin{Proof}
  Let $u,v$ be integers satisfying the conditions of Lemma~\ref{lemma_gcd_special}. Due to Lemma~\ref{lemma_lb_msirep} it suffices to
  prove that the optimal solution $(q',w_1',w_2')$ of LP~(\ref{lp_minsum}) 
  has a target value of at least $2ab$. Since $(v,u)$ is a losing vector and $(a,0)$, $(0,b)$ are winning vectors we have
  $uw_2'+vw_1'\le q'-1$, $aw_1'\ge q'$, and $bw_2'\ge q'$. Multiplying the first inequality by $ab$ yields
  $$
    abuw_2'+abvw_1'\le abq'-ab.
  $$
  Adding $bv$ times the second inequality and $au$ times the third inequality yields
  $$
    abuw_2'+abvw_1'\ge \underset{=ab-1}{\underbrace{(au+bv)}}q'.
  $$
  Thus we conclude $abq'-q'\le abq'-ab$, which is equivalent to $q'\ge ab$. 
  Next we deduce $w_1'\ge b$ and $w_2'\ge a$ from $aw_1'\ge q'\ge ab$ and $bw_2'\ge q'\ge ab$. Thus we have $aw_1'+bw_2'\ge 2ab$.
\end{Proof}

\begin{Corollary}
  \label{cor_subgame_lb}
  Let $\chi$ be a weighted voting game with equivalence classes $N_1,\dots,N_t$, \highlight{$1\le i_1<i_2\le t$ be two indices, and $N_1'\subset N_{i_1}$,
  $N_2'\subset N_{i_2}$ be two subsets. Consider two coprime integers $b>a\ge 1$, such that $|N_1'|=a$ and $|N_2'|=b$.} If the restriction
  of $\chi$ to $N_1'\cup N_2'$ is equivalent to $\chi_{a,b}$, then we have $q\ge ab$, $w_{i_1}\ge b$, and $w_{i_2}\ge a$ for the optimal
  solution $(q,w_1,\dots,w_t)$ of the linear program~(\ref{lp_minsum}).
\end{Corollary}

\subsection{Weighted games without a {\mirep} for small $t$}
\label{subsec_3_1}

In order to construct a weighted game without a {\mirep} for $t=3$ equivalence classes of voters we choose two coprime
integers $b>a\ge 1$ and an integer $c$ satisfying
\begin{enumerate}
  \item[(1)] $ab-c\notin\langle a,b\rangle$,
  \item[(2)] $ab-2c+1\in\langle a,b\rangle$,
  \item[(3)] $ab\ge 2c-1$, and
  \item[(4)] $c\ge b+\highlight{2}$.
\end{enumerate}
With this we consider the weighted game
$$
  \chi_{a,b,c}=[ab\highlight{;}c-\frac{1}{2},c-\frac{1}{2},\underset{a}{\underbrace{b,\dots,b}},\underset{b}{\underbrace{a,\dots,a}}],
$$
i.e.\ in the notation of a complete simple game the cardinality vector is given by $\cardv=(2,a,b)$.

These technically looking constraints can be interpreted as follows. Due to $b\ge a+1\ge 1$ and $c\ge b+\highlight{2}$ the
assignment $q=ab$, $w_1=a$, $w_2=b$, $w_3=c\highlight{-\frac{1}{2}}$ is a feasible solution of inequality system~(\ref{ie_feas2}) for $\chi_{a,b,c}$.
We remark $w_1=c-\frac{1}{2}\notin\mathbb{N}$. Constraint~(1) requires that every (\highlight{\smax}) losing vector $\cv{l}=(1,l_2,l_3)$ has a weight of
at most $q-\frac{3}{2}$, while constraints (2) and (3) ensure that there exists a winning vector $\highlight{\cv{m}_t}=(2,\highlight{t_2},\highlight{t_3})$
whose weight equals exactly the quota $q$, i.e.\ $2c-1+\highlight{t_2}b+\highlight{t_3}a=ab$.

\begin{Lemma}
  \label{lemma_ex_1}
  The sum of weights of a {\msirep} of $\chi_{a,b,c}$ is at least $2c-1+2ab$.
\end{Lemma}
\begin{Proof}
  Let $(q',w_1',w_2',w_3')$ be the optimal solution of the linear program minimizing the sum of weights. From Corollary~\ref{cor_subgame_lb} we
  conclude $q'\ge ab$, $w_2'\ge b$, and $w_3'\ge a$. Since \highlight{the, above defined, vector} $\cv{m}_t=(2,\highlight{t_2},\highlight{t_3})$,
  \highlight{is winning} we have $2w_1'+\highlight{t_2}w_2'+\highlight{t_3}w_3'\ge q'\ge ab$. \highlight{Using $b\highlight{t_2}+a\highlight{t_3}=ab-2c+1$}
  we conclude
  \begin{eqnarray*}
   && 2w_1'+aw_2'+bw_3'= 2w_1'+\highlight{t_2}w_2'+\highlight{t_3}w_3'+(a-\highlight{t_2})w_2'+(b-\highlight{t_3})w_3'\\
   &\ge& ab+(a-\highlight{t_2})b+(b-\highlight{t_3})a = 2ab+ab-\highlight{t_2}b\highlight{-}\highlight{t_3}a=2c-1+2ab
  \end{eqnarray*}
  \highlight{and finally apply Lemma~\ref{lemma_lb_msirep}.}
\end{Proof}

Next we show that $\tau_1=(ab,c,c-1,b,\dots, b,a,\dots a)$ and $\tau_2=(ab,c-1,c,b,\dots,b,a,\dots,a)$ are
{\msirep s} of $\chi_{a,b,c}$. Due to Lemma~\ref{lemma_ex_1} it remains to show that both vectors are {\irep s}. Coalitions of type
$(0,m_2,m_3)$ or $(2,m_2,m_3)$ have the same weight according to all three different weight vectors (including
$\tau_0=(c-\frac{1}{2},c-\frac{1}{2},b,\dots,b,a,\dots,a)$). Now let $(1,m_2,m_3)$ be a winning vector. From the definition of the
game $\chi_{a,b,c}$, i.e.\ $\tau_0$, we conclude $c-\frac{1}{2}+m_2b+m_3a\ge ab$, which can be slightly sharpened to $c-1+m_2b+m_3a\ge ab$.
Thus for both weightings $\tau_1$ and $\tau_2$ the lightest coalition, i.e.\ the one with minimal weight, of type $(1,m_2,m_3)$ has at
least a weight of $c-1+m_2b+m_3a\ge ab$. Now
let $(1,l_2,l_3)$ be a losing coalition. Since $ab-c\notin\langle a.b\rangle$ we have $c-\frac{1}{2}+bl_2+al_3\le ab-\frac{3}{2}$. Thus
for both weightings $\tau_1$ and $\tau_2$ the heaviest coalition of type $(1,l_2,l_3)$ has a weight of at most a weight
of $c+m_2b+m_3a\le ab-1$.

The final conclusion is that $\chi_{a,b,c}$ cannot admit a {\mirep} since it has at least two {\msirep s}. An example where
the requested technical conditions on $a$, $b$, and $c$ are satisfied is given by $c=12$, $b=7$, $a=5$.

Instead of using the general Lemma~\ref{lemma_ex_1} we can consider the example directly. The game $\chi_{5,7,12}$ is a
complete simple game with $t=3$ types of voters whose vector/matrix notation is given by $\cardv=(2,5,7)$ and 
$\mathcal{M}=\begin{pmatrix}
2&2&1&1&1&0&0&0\\
1&0&2&1&0&5&3&0\\
1&3&2&4&5&0&3&7
\end{pmatrix}^T$. The matrix of the \highlight{\smax} losing vectors is given by 
$\begin{pmatrix}
2&2&1&1&0&0\\
1&0&3&1&4&2\\
0&2&0&3&1&4
\end{pmatrix}^T$. Solving the LP~(\ref{lp_minsum})
yields the optimal
solution $(35,11.5,7,5)$. Thus we conclude from Lemma~\ref{lemma_lb_msirep} that the sum of the weights in a {\msirep} is at least $93$. Now we can
easily check that both $(35,12,11,7,7,7,7,7,5,5,5,5,5,5,5)$ and $(35,11,12,7,7,7,7,7,5,5,5,5,5,5,5)$ are {\irep s} of $\chi_{5,7,12}$ attaining
this lower bound. We remark that the stated representations arise by a swap of weights within the first equivalence class.

\medskip

For $t=4$ equivalence classes and for the situation of {\ireppts}, we apply a similar idea and consider a game $\chi$ with cardinality
vector $\cardv=(1,1,a,b)$
containing $\chi_{a,b}$ as a subgame. The rough idea is to choose half-\highlight{integer} weights $w_1,w_2\in\mathbb{N}+\frac{1}{2}$ such that $(ab,w_1,w_2,b,a)$
is an optimal solution of LP~(\ref{lp_minsum}) 
while $\lfloor w_1\rfloor$ or $\lfloor w_2\rfloor$ can be attained in (different) {\irep s}. Similarly as in the example above, sufficient technical conditions can be formulated
using the membership or non-membership of certain values in $\langle a,b\rangle$. We refrain from explicitly formulating the details and instead give an example.
We choose $b=11$ and $a=7$, which satisfy $\{52,59\}\cap\langle 7,11\rangle=\emptyset$ and $52+59-7\cdot 11+1\in\langle 7,11\rangle$. The game $\chi$ now
is uniquely chosen by stating its matrix of \highlight{\smin} winning vectors:
$\mathcal{M}=\left(\begin{array}{rrrrrrrrrrrrrrr}
1&1&1&1&1&1&0&0&0&0&0&0&0&0&0\\
1&1&0&0&0&0&1&1&1&1&0&0&0&0&0\\
2&0&5&3&1&0&5&3&1&0&7&6&4&2&0\\
2&5&0&3&6&8&1&4&7&9&0&2&5&8&11
\end{array}\right)^T\!$. The matrix of the \highlight{\smax} losing vectors is then given by
$\left(\begin{array}{rrrrrrrrrrrrr}
1&1&1&1&1&0&0&0&0&0&0&0&0\\
1&1&0&0&0&1&1&1&1&0&0&0&0\\
3&1&4&2&0&5&4&2&0&6&5&3&1\\
0&3&1&4&7&0&2&5&8&1&3&6&9
\end{array}\right)^T\!$.
Solving the LP~(\ref{lp_minsum})
yields the optimal solution $(77,24.5,17.5,11,7)$.
The sum of weights of a {\msirep} is at least $196$. By checking that both $(77,25,17,11,11,11,11,11,11,11,$
$7,7,7,7,7,7,7,7,7,7,7)$ and $(77,24,18,11,11,11,11,11,11,11,7,7,7,7,7,7,7,7,7,7,7)$ are {\irep s} of $\chi$ we conclude that they are indeed {\msireppts}. Thus $\chi$ does not admit
a {\mireppt}.

\subsection{Weighted games with more than two {\msirep s}}
\label{subsec_3_2}
It would be nice to have an example of a weighted game with more than two minimum sum representations preserving types.
Before we give a construction that works, we briefly remark that not every reasonable system of constraints for the
representability of some expressions needs to have a solution, so that not all construction ideas lead to success. Our first idea
was to choose $t=4$, two coprime integers $b>a\ge 1$ and $\cardv=(1,1,a,b)$. If there existed integers
$0<l_2<l_1-1<ab$ such that $l_1,l_1+1,l_2,l_2+1\notin\langle a,b\rangle$ but $ab-l_1-l_2+2\in\langle a,b\rangle$,
then we could check that $(ab,ab-l_2-2+x,ab-l_1-2+(2-x),b,\dots,b,a,\dots,a)$ is a {\msirep} for $x\in\{0,1,2\}$. Unfortunately
the existence of such integers $l_1,l_2$ would contradict Popoviciu's theorem \highlight{counting the number of representations, see e.g.\ 
\cite{1114.52013}}. \highlight{To be more precisely,} \ $ab-k\in\langle a,b\rangle$ implies $k\notin\langle a,b\rangle$ \highlight{for all
$k\in\mathbb{N}$ with $a,b\!\not| k$}.

For $t=5$ we have another construction which works:
\begin{Proposition}
   Let $b>a\ge 1$ be two coprime positive integers. Suppose we have integers $l_1<l_2<l_3$ fulfilling
   \begin{itemize}
     \item[(1)] $b<\widehat{w_i}=ab-l_i-1$,
                       $l_i\notin\langle a,b\rangle$ 
                for $1\le i\le 3$,
     \item[(2)] $l_1+l_2+l_3-l_i-ab\notin\langle a,b\rangle$ for $1\le i\le 3$, and
     \item[(3)] $0<l_1+l_2+l_3+1-2ab<ab$, $l_1+l_2+l_3+1-2ab\in\langle a,b\rangle$.
   \end{itemize}
   With this the weighted game $\chi=[ab\highlight{;}\widehat{w_1},\widehat{w_2}+1,\widehat{w_3}+1,b,\dots,b,a,\dots,a]$,
   where $\cardv=(1,1,1,a,b)$, has the following three {\msireppts}:
   \begin{itemize}
    \item $\tau_1=(ab,\widehat{w_1},\widehat{w_2}+1,\widehat{w_3}+1,b,\dots,b,a,\dots,a)$
    \item $\tau_2=(ab,\widehat{w_1}+1,\widehat{w_2},\widehat{w_3}+1,b,\dots,b,a,\dots,a)$
    \item $\tau_3=(ab,\widehat{w_1}+1,\widehat{w_2}+1,\widehat{w_3},b,\dots,b,a,\dots,a)$
   \end{itemize}
\end{Proposition}
\begin{Proof}
  \highlight{Let} $(q',w_1',w_2',w_3',w_4',w_5')$ be the optimal solution of the linear
  program minimizing the sum of weights. From Corollary~\ref{cor_subgame_lb} we conclude $q'\ge ab$, $w_4'\ge b$, and $w_5'\ge a$. \highlight{Since}
  $l_1+l_2+l_3+1-2ab\in\langle a,b\rangle$, see constraint~(3), there exist integers $u,v$ such that $(1,1,1,u,v)$ is a \highlight{(\smin)}
  winning vector and $ub+va=l_1+l_2+l_3+1-2ab$. Thus we have $w_1'+w_2'+w_3'+uw_4'+vw_5'\ge q'$. Inserting this into the sum of
  all weights yields
  \begin{eqnarray*}
    && w_1'+w_2'+w_3'+aw_4'+bw_5' \ge q'+(a-u)w_4'+(b-v)w_5'\\
    &\ge& ab+2ab-(ub+va)= 2ab+\widehat{w}_1+\widehat{w}_2+\widehat{w}_3+2,
  \end{eqnarray*}
  i.e.\ each {\msirep} has a weight of at least $2ab+\widehat{w}_1+\widehat{w}_2+\widehat{w}_3+2$ \highlight{due Lemma~\ref{lemma_lb_msirep}}.
  
  The next step is to prove that the three stated weightings represent the game $\chi$. So for each vector $(m_1,m_2,m_3,m_4,m_5)$
  we have to prove that its weight is less then $ab$ or at least $ab$ in all three different weightings simultaneously. This can be
  \highlight{easily verified} for the cases where $m_1+m_2+m_3\in\{0,3\}$.
  
  For $m_1+m_2+m_3=1$ let $(1,0,0,m_4,m_5)$ be a losing vector, i.e.\ we have $\widehat{w}_1+m_4b+m_5a=
  ab-l_1-1+m_4b+m_5a\le ab-1$. In $\tau_2$ and $\tau_3$ the weight of the first player is increased by one, so we need $ab-l_1-1+m_4b+m_5a\le ab-2$,
  which is indeed true since $ab-l_1-1+m_4b+m_5a= ab-1$ is equivalent to $l_1\in\langle a,b\rangle$, i.e.\ it contradicts constraint~(1). Increasing
  the weight of player one in a winning coalition does not \highlight{affect its status}. Now we consider a winning vector $(0,1,0,m_4,m_5)$, i.e.\ we have
  $\widehat{w}_2+1+m_4b+m_5a=ab-l_2+m_4b+m_5a\ge ab$. In $\tau_2$ the weight of the second player is decreased by one, so that we need
  $ab-l_2+m_4b+m_5a\ge ab+1$, which is true since $ab-l_2+m_4b+m_5a= ab$ is equivalent to $l_2\in\langle a,b,\rangle$. Due to symmetry we conclude
  that there are no contradictions for $m_3=1$.

  \highlight{
  For $m_1+m_2+m_3=2$ let $(1,1,0,m_4,m_5)$ be a losing vector, i.e.\ we have $\widehat{w}_1+\widehat{w}_2+1+m_4b+m_5a=
  2ab-l_1-l_2-1+m_4b+m_5a\le ab-1$. Due to $l_1+l_2\notin\langle a,b\rangle$ the vector has a weight of at most $ab-2$ using the weights from $\tau_1$. Thus the
  vector remains losing in $\tau_2$ and $\tau_3$. Increasing the weight of player one in a winning coalition does not \highlight{affect its status}.
  A symmetric argument
  applies for vectors of type $(1,0,1,m_4,m_5)$. Now let $(0,1,1,m_4,m_5)$ be a winning vector, i.e.\ we have $\widehat{w}_2+\widehat{w}_3+2+m_4b+m_5a=
  2ab-l_1-l_2+m_4b+m_5a\ge ab$. Due to $l_1+l_2\notin\langle a,b\rangle$ the vector has a weight of at least $ab+1$ using the weights from $\tau_1$.
  Thus the vector remains winning in $\tau_2$ and $\tau_3$. Decreasing the weight of either player two or player three by one does not cause any problems
  for a losing coalition.
  }
  
  Thus we have exhaustively checked that $\tau_1$, $\tau_2$, and $\tau_3$ are {\msireppts}.
\end{Proof}

An example where the requirements of the previous proposition are fulfilled is given by $a=17$, $b=13$, $l_1=157$, $l_2=161$, $l_3=174$, $\widehat{w_1}=63$, $\widehat{w_2}=59$, and $\widehat{w_3}=46$.  A smaller example is given by $a=13$, $b=11$, $l_1=93$, $l_2=97$, $l_3=106$, $\widehat{w_1}=49$,
$\widehat{w_2}=45$, and $\widehat{w_3}=36$. Furthermore we have the following straightforward generalization:

\begin{Proposition}
   Let $b>a\ge 1$ be two coprime integers with $a>b$ and $t$ be an integer with $t\ge 2$.
   Suppose we have integers $l_1<l_2<\dots<l_t$ fulfilling
   \begin{itemize}
     \item[(1)] $b<\widehat{w_i}=ab-l_i-1$, $l_i\notin\langle a,b\rangle$ for $1\le i\le 3$,
     \item[(2)] $\sum\limits_{j=1}^z l_{i_j}-(z-1)ab\notin\langle a,b\rangle$  for all $2\le z< t$ and all subsets
                $\left\{i_1,\dots,i_z\right\}\subseteq\{1,\dots,t\}$ of cardinality $z$, and
     \item[(3)] $0<\sum\limits_{j=1}^t l_j\,+1-(t-1)ab<ab$, $\sum\limits_{j=1}^t l_j\,+1-(t-1)ab\in\langle a,b\rangle$.
   \end{itemize}
   With this the weighted game $\chi=[ab\highlight{;}\widehat{w_1},\widehat{w_2}+1,\dots,\widehat{w_t}+1,b,\dots,b,a,\dots,a]$,
   where $\cardv=(\underset{t}{\underbrace{1,\dots,1}},a,b)$,
   has the following $t$ {\msireppts}:
   \begin{itemize}
    \item $(ab,\widehat{w_1},\widehat{w_2}+1,\dots,\widehat{w_t}+1,b,\dots,b,a,\dots,a)$
    \item $(ab,\widehat{w_1}+1,\widehat{w_2},\widehat{w_3}+1\dots,\widehat{w_t}+1,b,\dots,b,a,\dots,a)$\\[0.5mm]
          \hspace*{2mm}$\vdots$
    \item $(ab,\widehat{w_1}+1,\dots,\widehat{w_{t-1}}+1,\widehat{w_t},b,\dots,b,a,\dots,a)$
   \end{itemize}
\end{Proposition}

For $t\ge 4$ we have the following examples:\\[-5mm]
\begin{center}
  \begin{tabular}{rrrl}
    \hline
    $\mathbf{t}$ & $\mathbf{a}$ & $\mathbf{b}$ & $\mathbf{l_1,\dots,l_t}$\\
    4 & 19 & 11 & 141, 157, 160, 179\\
    5 & 19 & 17 & 249, 251, 253, 268, 287\\
    6 & 29 & 17 & 389, 396, 401, 418, 430, 447\\
    7 & 31 & 29 & 746, 750, 752, 777, 779, 808, 810\\
    8 & 37 & 29 & 883, 891, 920, 941, 949, 970, 978, 1007\\
    9 & 41 & 31 & 1086, 1100, 1106, 1117, 1127, 1137, 1158, 1168, 1199\\
   10 & 43 & 41 & 1513, 1550, 1552, 1554, 1593, 1595, 1597, 1636, 1638, 1679\\
    \hline
  \end{tabular}
\end{center}


\highlight{
\begin{Conjecture}
  Each weighted game with $t$ equivalence classes of voters can have at most $t-2$ different minimum sum representations
  preserving types.
\end{Conjecture}}

\subsection{Possible weights of {\mirep s}}
\label{subsec_3_3}
Instead of asking which classes of weighted games admit a {\mirep} or a {\mireppt} one can ask which weights are possible
in a {\mirep}. The following theorem and remarks resolve this question for two different weights almost completely. The stated
lower bounds on the number of necessary voters $n$ might be improved.

\begin{Theorem}
  \label{thm_min_rep_two_weights}
  For two coprime integers $b>a\ge 1$ 
  the weighted game $\chi=[q=ab\highlight{;}\underset{n_1}{\underbrace{b,\dots,b}},\underset{n_2}{\underbrace{a,\dots,a}}]$,
  where $n_1\ge \highlight{a}$ and $n_2\ge \highlight{b}$, is in {\mirep}.
\end{Theorem}
\begin{Proof}
  Let $(q',b_1',\dots,b_{n_1}',a_1',\dots,a_{n_2}')$ be an arbitrary {\irep} of $\chi$, where we assume $b_1'\ge\dots\ge b_{n_1}'$ and
  $a_1'\ge\dots\ge a_{n_2}'\ge 0$ w.l.o.g. From Isbell's desirability relation we conclude $b_{n_1}'>a_1'$. \highlight{By} 
  Corollary~\ref{cor_subgame_lb} every {\irep} of $\chi$ has a sum of weights of at least $n_1b+n_2a$ \highlight{and $q'\ge ab$, 
  $\sum_{i=1}^{n_1} b_i'\ge n_1b$, $\sum_{i=1}^{n_2} a_i'\ge n_2a$. It suffices to show $b_{n_1}'\ge b$ and $a_{n_2}'\ge a$.}
  
  \highlight{If $a_{n_2}'<a$ we can assume $a_{n_2}'=a-1$, since convex combinations of feasible weightings are feasible. By averaging
  the weights $a_1',\dots a_{n_2-1}'$ and $b_1',\dots,b_{n_1}'$ we obtain the (feasible, possibly \highlight{non-integer}) weighting 
  $(q',b+t,\dots,b+t,a+s,\dots,a+s,a-1)$, where $s\in\mathbb{Q}_{>0}$, $t\in\mathbb{Q}_{\ge 0}$.}
 
  Let $0\le u\le a-1$ and $0\le v\le b-1$ be two integers with $ub+va=q-1=ab-1$. \highlight{Rearranging yields $u=a-\frac{1+av}{b}$ so that
  $b$ divides $1+av$ and we have $b<av$.} Since $(u,v)$ is a losing vector and $(a,0)$, $(0,b)$ are winning vectors we have
  \highlight{\begin{equation}
    \label{ie_possible_weights}
    ub+ut+av+vs\le q'-1,\quad
    ab+at\ge q',\,\,\text{ and }\,\,
    ab+bs-s-1\ge q'.
  \end{equation}}
  Multiplying the first inequality by $ab$ yields
  \highlight{
  $$
    ab^2u+abut+a^2bv+abvs\le abq'-ab
  $$}
  and $bu$ times the second inequality plus $av$ times the third inequality yields
  \highlight{
  $$
    ab^2u+abtu+
    a^2bv+absv-avs-av
    \ge abq'-q'
  $$}
  \highlight{Combining the last two inequalities yields
  \begin{equation}
    q'\ge ab+(s+1)av. 
  \end{equation}
  We already know $ab+bs-s-1\ge q'$ and conclude $(b-1)s-1\ge (s+1)av$. Inserting $b<av$ yields the contradiction $-s-1>b$. 
  Thus $a_{n_2}'\ge a$.
  }
  
  \medskip
  
  \highlight{
  If $b_{n_1}'<b$ we can assume $b_{n_1}'=b-1$ and consider the weighting 
  $(q',b+t,\dots,b+t,b-1,a+s,\dots,a+s)$, with $t\in\mathbb{Q}_{>0}$, $s\in\mathbb{Q}_{\ge 0}$. Let again $0\le u\le a-1$ and
  $0\le v\le b-1$ be two integers with $ub+va=q-1=ab-1$. Here we have $a<ub$. An analogous calculation as before yields
  $$
    q'\ge ab+(t+1)ub
  $$
  and
  $$
    q'\le ab+(a-1)t-1.
  $$
  Combining these two inequalities yields $(a-1)t-1\ge(t+1)ub$. Inserting $a<ub$ ends up in the contradiction $-t-1>a$.
  }
\end{Proof}

In the following remark we want to emphasize that most of the requirements of Theorem~\ref{thm_min_rep_two_weights} are necessary:

\begin{Remark} $\,$\\[-3mm]
  \begin{itemize}
   \item[(1)] If $r=\gcd(a,b)>1$, then $\left[\frac{q}{r},\frac{b}{r},\dots,\frac{b}{r},\frac{a}{r},\dots,\frac{a}{r}\right]$ is a smaller representation
              for the same game.
   \item[(2)] If $b=a$ then there is only one type of voters with minimum representation $[q'\highlight{;}1,\dots,1]$ for a suitable quota $q'$. If $b<a$ then the voters
              of type~$2$ would be more powerful than the voters of type~$1$, which is not possible by definition.
   \item[(3)] If $a=0$ and $b>1$ then $\left[\left\lceil\frac{q}{b}\right\rceil,1,\dots,1,0,\dots,0\right]$ is a smaller representation for the same game,
              which is indeed the {\mirep}.
   \item[(4)] The lower bounds on $n_1$ and $n_2$ can be improved, e.g.\ based on the knowledge of $u$ and $v$.
  \end{itemize}
\end{Remark}

There is a generalization to weighted games with more than two types of voters:
\begin{Theorem}
  Let $a_1,\dots,a_t$ be integers such that $a_1>a_2>\dots>a_t>0$ and for each $1\le i\le t$ there is an index $1\le j\le t$ with $\gcd(a_i,a_j)=1$.
  The weighted game $$\chi=[q=lcm(a_1,\dots,a_t)\highlight{;}\underset{n_1}{\underbrace{a_1,\dots,a_1}},\dots,\underset{n_t}{\underbrace{a_t,\dots,a_t}}],$$
  where $n_i\ge \highlight{lcm(a_1,\dots,a_t)}/a_i$ for all $1\le i\le t$, is in {\mirep}.
\end{Theorem}
\begin{Proof}
  For an arbitrary {\irep} of $\chi$ let $(q',a_1',\dots,a_t')$ be the averaged representation with equal (possibly \highlight{non-integer}) weights 
  within each equivalence class of voters. For each index $1\le i\le t$ choose a suitable index $j\neq i$ such that $\gcd(a_i,a_j)=1$.
  Let $z=\frac{q}{a_ia_j}\in\mathbb{N}$ and
  $0\le u\le a_j-1$, $0\le v\le a_i-1$ be two integers such that $ua_i+va_j=a_ia_j-1$.
  With this we have
  \begin{eqnarray*}
    -\Big(\left(z-1\right)a_j+u\Big)a_i'-va_j'+q' &\ge& 1\\
    (z-1)a_ja_i'+a_ia_j'-q' &\ge& 0\\
    za_ja_i'-q' &\ge& 0.
  \end{eqnarray*}
  Combining these inequalities with the vectors $(a_i,v,a_i-v)$, $(a_j,a_j-u,u)$, and\\
  $(za_ia_j,za_jv,(z-1)a_j(a_i-v)+a_iu)$ as multipliers yields $a_i'\ge a_i$, $a_j'\ge a_j$, and $q'\ge q$.
  
  We can treat the case of different weights within equivalence classes of voters analogously to the proof of the previous theorem.
\end{Proof}

In the next theorem we pay for less restrictive conditions on the weights $a_i$ by a rather large bound on the number of voters $n$.
To this end we generalize Lemma~\ref{lemma_gcd_special} for more than two integers:
\begin{Lemma}
  \label{lemma_xgcd_general}
  Let $a_1,\dots,a_t$ be positive integers with $t\ge 2$ and $\gcd(a_1,\dots,a_t)=g$. There exist $t$ integers $u_i$ with
  $\sum\limits_{i=1}^tu_ia_i=\prod_{i=1}^t a_i-g$ and $0\le u_i\le \prod\limits_{j=1,j\neq i}^t a_j\,-1$ for all $1\le i\le t$.
\end{Lemma}
\begin{Proof}
  We prove by induction on $t$. For $t=2$ we apply Lemma~\ref{lemma_gcd_special} for the two integers $\frac{a_1}{g}$, $\frac{a_2}{g}$. For $t>2$
  let $u_i'$ be integers with $\sum\limits_{i=1}^{t-1}u_i'a_i=\prod_{i=1}^{t-1} a_i-g'=:k$ and $0\le u'_i\le \prod\limits_{j=1,j\neq i}^{t-1} a_j\,-1$
  for all $1\le i\le t-1$, where $g'=\gcd(a_1,\dots,a_{t-1}')$. We remark $\gcd(k,a_t)=g$ and apply Lemma~\ref{lemma_xgcd_general} for $t=2$.
\end{Proof}

\begin{Theorem}
  For integers $a_1>a_2>\dots>a_t>0$ with $\gcd(a_1,\dots,a_t)=1$ the weighted game
  $$\chi=[q=\prod\limits_{j=1}^t a_j\highlight{;}\underset{n_1}{\underbrace{a_1,\dots,a_1}},\dots,\underset{n_t}{\underbrace{a_t,\dots,a_t}}],$$
  where $n_i\ge 2\prod\limits_{j=1,j\neq i}^t a_j$ for all $1\le i\le t$, is in {\mirep}.
\end{Theorem}
\begin{Proof}
  Due to Lemma~\ref{lemma_xgcd_general} there are $t$ integers $u_i$ with $\sum\limits_{i=1}^tu_ia_i=q-1$
  and $0\le u_i\le \prod\limits_{j=1,j\neq i}^t a_j\,-1$ for all $1\le i\le t$. For an arbitrary {\irep} 
  of $\chi$ let $(q',a_1',\dots,a_t')$ be the averaged representation with equal (possibly \highlight{non-integer}) weights 
  within each equivalence class of voters. With this the following inequalities have to be valid:
  $$
    \frac{q}{a_i}\cdot a_i'-q' \ge 0\text{ for all } 1\le i\le t\quad\text{and}\quad
    -\sum_{i=1}^t u_ia_i' +q' \ge 1.
  $$
  Summing up $a_iu_i$ times the $i$th inequality plus $q$ times the last inequality yields $q'\ge q$. Inserting this into the
  $i$th inequality gives $a_i'\ge a_i$ for all $1\le i\le t$.
  
  We can treat the case of different weights within equivalence classes of voters analogously to the proofs of the previous theorems.
\end{Proof}

\highlight{The} condition $\gcd(a_1,\dots,a_t)=1$ is necessary. If we also want to use zero weights we can utilize the next lemma:

\begin{Lemma}
  The weighted game $[q\highlight{;}\underset{n_1}{\underbrace{a_1,\dots,a_1}},\dots,\underset{n_t}{\underbrace{a_t,\dots,a_t}}]$, where $a_i>0$ for all
  $1\le i\le t$, is in {\mirep} \emph{if and only if} the weighted game
  $[q\highlight{;}\underset{n_1}{\underbrace{a_1,\dots,a_1}},\dots,\underset{n_t}{\underbrace{a_t,\dots,a_t}},\underset{n_{t+1}}{\underbrace{0,\dots,0}}]$ is in
  {\mirep}.
\end{Lemma}

\section{Weighted games with {\mirep s} for two types of voters}
\label{sec_4}

\noindent
As we have remarked in the introduction, all complete simple games with just one type of voters, $t=1$, are weighted and admit a {\mirep} with
all weights being equal to~$1$. These games are called ``symmetric``, ``anonymous`` or ``$k$-out-of-$n$-games`` in the literature. In the previous section we have
constructed weighted games with $t=3$ equivalence classes without a {\mirep}. So the central question of this section (and the paper) is:
what happens for two types of voters? We first state the main result.

\begin{Theorem}
  \label{main_thm}
  Each weighted game with two types of voters admits a {\mirep} $(q,\overset{n_1}{\overbrace{w_1,\dots,w_1}},\overset{n_2}{\overbrace{w_2,\dots,w_2}})$,
  with $1\le w_1\le \max(n_1+1,n_2)$, $0\le w_2\le \max(n_1,n_2-1)$, and
  $1\le q\le (n_1+n_2-1)\cdot\max(n_1+1,n_2)$. For $r\ge 2$ \highlight{\smin} winning vectors
  the bounds of minimum weights can be sharpened to $1\le w_1\le n_2$, $1\le w_2\le n_1$, and $w_2+1\le q\le 2n_1n_2$.
\end{Theorem}

In Section~\ref{sec_5} we will use the bounds on $q$, $w_1$, and $w_2$ to give an upper bound on the number of weighted voting games
with two types of voters. As preliminary work we prove Theorem~\ref{main_thm} for some special cases of weighted games with two types of voters by a
direct argumentation on the possible {\irep s} in Subsection~\ref{subsec_4_1}. The proof strategy for the remaining part is more involved. 
We study linear minimization problems subject to the constraints in (\ref{ie_feas2}). It turns out that each target function without negative
coefficients admits an optimal integer solution. By additionally using some structure result on the set of inequalities, which attain equality,
called \textit{tight} later on, we deduce the existence of a {\mirep}.

\subsection{Proof of the main theorem for $r=1$ \highlight{shift-minimal winning vectors}}
\label{subsec_4_1}

\begin{Theorem}
  \label{thm_min_rep_2_1}
  For a weighted game $\chi$ with two types of voters, a cardinality vector $\cardv=(n_1,n_2)$ and a unique minimal winning coalition
  $\widetilde{m}=(m_1,m_2)$, i.e.\ $t=2$ and $r=1$, there exists a {\mirep}.
\end{Theorem}
\begin{Proof}
  Since the game is weighted there are some restrictions on the parameters $m_1$, $m_2$ beyond those from (\ref{compact_ilp_2_1}), i.e.\
  $1\le n_1\le n-1$, $n_1+n_2=n$, $1\le m_{1} \le n_1$, and $0\le m_{2}\le n_2-1$.
  \highlight{First we exclude} the cases where $1\le m_1\le n_1-1$ and $2\le m_2\le n_2-2$\highlight{. A}ssume that $(q,w_1,w_2)$ is a feasible
  solution of (\ref{ie_feas2}). Since $(m_1-1,m_2+2)$, $(m_1+1,m_2-2)$ are losing vectors and $(m_1,m_2)$ is a winning vector we have
  \begin{eqnarray*}
    (m_1-1)w_1+(m_2+2)w_2\le q-1\le m_1w_1+m_2w_2-1,\\
    (m_1+1)w_1+(m_2-2)\highlight{w_2}\le q-1\le m_1w_1+m_2w_2-1,
  \end{eqnarray*}
  from which we conclude $2w_2\le w_1-1\le 2w_2-2$; a contradiction. It will turn out that $\chi$ is weighted in the remaining cases, i.e.\
  for $m_1=n_1$ or $m_2\in\{0,1,n_2-1\}$.

  Let $(q,a_1,\dots,a_{n_1},b_1,\dots,b_{n_2})$ be an arbitrary {\irep} with $a_1\ge\dots\ge a_{n_1}$ and $b_1\ge \dots\ge b_{n_2}$ of
  the game $\chi$.  Due to Isbell's desirability relation we have $a_{n_1}\ge b_1+1$.

  \begin{itemize}
   \item $1\le m_1\le n_1-1$, $m_2=0$:\\
         We can easily check that $a_1=\dots=a_{n_1}=1$, $b_1=\dots=b_{n_2}=0$, $q=m_1$ is an {\irep} of $\chi$. Since we
         have $a_i\ge 1$ and $b_j\ge 0$ it is also a {\mirep}.
   \item $1\le m_1\le n_1-1$, $m_2=1$:\\
         Since $(m_1,0)$ and $(m_1-1,n_2)$ are losing vectors, we can conclude $b_i\ge 1$ for all $1\le i\le n_2$ and $a_i\ge n_2$ for all $1\le i\le n_1$.
         We can easily check that $a_1=\dots=a_{n_1}=n_2$, $b_1=\dots=b_{n_2}=1$, $q=m_1n_2+1$ is an {\irep} of $\chi$ and thus is 
         indeed a {\mirep}.
   \item $1\le m_1\le n_1-1$, $m_2=n_2-1$:\\
         Since the cases $m_2\in\{0,1\}$ were dealt previously, we assume $m_2\ge 2$.
         
         For $m_1+n_2-1\le n_1$ the vector $(m_1+n_2-2,0)$ is losing. Comparing the weights of its
         corresponding coalitions with those from the \highlight{\smin} winning vector and inserting $a_i\ge b_j+1$ yields $b_j\ge n_2-1$ and $a_i\ge n_2$. We can
         easily check that $a_1=\dots=a_{n_1}=n_2$, $b_1=\dots=b_{n_2}=n_2-1$, $q=m_1n_2+(n_2-1)^2$ is an {\irep} of $\chi$.
         
         \medskip
         
         \noindent
         For $m_1+n_2-1>n_1$ we compare the losing vector $(n_1,m_1+n_2-2-n_1)$ with the \highlight{\smin} winning vector
         and insert $a_i\ge b_j+1$ to deduce $b_j\ge n_1-m_1+1$ and $a_i\ge n_1-m_1+2$. We can easily check that $a_1=\dots=a_{n_1}=n_1+2-m_1$,
         $b_1=\dots=b_{n_2}=n_1+1-m_1$, $q=(m_1+n_2)(n_1+1-m_1)+2m_1-n_1-1$ is an {\irep} of $\chi$.
   \item $m_1=n_1$, $0\le m_2\le n_2-1$:\\
         If $m_2=0$ then $a_1=\dots=a_{n_1}=1$, $b_1=\dots=b_{n_2}=0$, $q=n_1$ is a {\mirep}. Otherwise we have 
         the losing vectors $(n_1,m_2-1)$ and $(n_1-1,n_2)$ from which we conclude $b_i\ge 1$ for all $1\le i\le n_2$ and 
         $a_i\ge n_2-m_2+1$, respectively. Thus we have $q\ge n_1(n_2-m_2+1)+m_2$. We can easily check that equality is possible, so that
         $a_i=n_2-m_2+1\ge 2$, $b_j=1$, $q=n_1(n_2-m_2+1)+m_2$ is a {\mirep}.
  \end{itemize}
\end{Proof}

Going over the cases of the proof of Theorem~\ref{thm_min_rep_2_1} we can check that all stated {\mirep s} satisfy $1\le w_1\le\max(n_1+1,n_2)$
and $0\le w_2\le \max(n_1,n_2-1)$ so that $1\le q\le (n_1+n_2-1)\cdot \max(n_1+1,n_2)$, as stated in Theorem~\ref{main_thm}.

To reduce the need for case differentiations in the remaining part we now completely handle the cases where null voters or dummies, i.e.\ voters
$i$ such that $\chi(U)=\chi(U\cup\{i\})$ for all subsets $U\subseteq N\backslash\{i\}$, occur.

\begin{Lemma}
  Weighted games with two types of voters, where one class consists of null voters, admit a {\mirep}.
\end{Lemma}
\begin{Proof}
  If, as usual, the equivalence classes of the game $\chi$ are given by $N_1$, $N_2$, then $N_2$ has to be the set of
  null voters. By the definition of a null voter each \highlight{\smin} winning vector $(m_1,m_2)$ has to satisfy $\highlight{m_2}=0$. Since
  \highlight{\smin} winning vectors are incomparable, we have $r=1$ and can apply Theorem~\ref{thm_min_rep_2_1}.
\end{Proof}

Also in general we can drop null voters from given games when determining {\mirep s} or {\msirep s} (preserving types or not).
\begin{Lemma}
  Let $\chi$ be a weighted game with $k$ null voters and $\chi'$ be the (weighted) \highlight{game} arising from $\chi$ by
  deleting the $k$ null voters. If $(q,w_1,\dots,w_{n-k})$
  is an {\irep} of $\chi'$, then $(q,w_1,\dots,w_{n-k},0,\dots,0)$ is an {\irep} of $\chi$.
\end{Lemma}

\subsection{Proof of the main theorem for $r>1$ \highlight{shift-minimal winning vectors}}
\label{subsec_4_2}

In the following we restrict our considerations to games without null voters and $r>1$ \highlight{shift-minimal winning vectors}. In
this case we can drop two constraints from inequality system~(\ref{ie_feas2}).

\begin{Lemma}
  \label{lemma_reduced_feas_lp}
  For a weighted game $\chi$ without null voters and with $t=2$, $r>1$ every vector $(q,w_1,w_2)$ is feasible for inequality system~(\ref{ie_feas2}) if and
  only if it satisfies
  \begin{equation}
     \cv{x}^T w\ge q\,\, \forall\, \cv{x}\in\mwvs,\,\,\cv{y}^Tw\le q-1\,\, \forall\, \cv{y}\in\mlvs.
     \label{ie_feas3}
  \end{equation}
\end{Lemma}
\begin{Proof}
  It remains to prove that the constraints from (\ref{ie_feas3}) imply $w_1\ge w_2+1$ and $w_2\ge 0$. 
  
  Let $(a,b)\in\mwvs$ with minimal $a$, i.e.\ all $(m_1,m_2)\in\mwvs$ satisfy $m_1\ge a$. If $a\ge 1$ and $b<n_2$, then
  $(a-1,n_2)\in\mlvs$. With this we conclude
  $$
    aw_1+bw_2\ge q\ge (a-1)w_1+n_2w_2+1,
  $$
  which is equivalent to $w_1\ge (n_2-b)w_2+1\ge w_2+1$.
  If $a=0$ or $b=n_2$ then let $(c,d)\in\mwvs$ with minimal $c>a$, i.e.\ for all $(m_1,m_2)\in\mwvs$ we either have $m_1=a$ or $m_1\ge c$.
  With this we have $(c-1,a+b-c)\in\mlvs$ and conclude
  $$
    cw_1+dw_2\ge q\ge (c-1)w_1+(a+b-c)w_2+1,
  $$
  which is equivalent to $w_1\ge (a+b-c-d)w_2+1\ge w_2+1$. Thus in both cases the constraints from (\ref{ie_feas3}) imply $w_1\ge w_2+1$.
  
  In order to deduce $w_2\ge 0$ we consider a winning vector
  $(m_1,m_2)$ with $m_2>0$, which must exist since $\chi$ does not contain null voters. Thus $(m_1,m_2-1)$ is a losing vector. Now let $(l_1,l_2)$ be a
  \highlight{\smax} losing vector with $(l_1,l_2)\succeq(m_1,m_2-1)$, i.e.\ we have $l_1\ge m_1$ and $l_1+l_2\ge m_1+m_2-1$. From $(l_1,l_2)\nsucceq (m_1,m_2)$ we
  conclude $l_1+l_2<m_1+m_2$, so that only $l_1+l_2=m_1+m_2-1$ is possible. With this we have
  $m_1w_1+\highlight{m_2}w_2\ge q\ge l_1w_1+l_2w_2+1$, which is equivalent to $(l_1-m_1)w_1+1\le (m_2-l_2)$. Inserting $w_1\ge w_2+1$ yields 
  $(l_1-m_1)w_2+1+(l_1-m_1)\le (m_2-l_2)w_2$, so that we have $w_2\ge 1+l_1-m_1\ge 1$.
\end{Proof}

Let us consider an example of inequality system (\ref{ie_feas3}) for the complete simple game $\chi$ uniquely characterized by $\cardv=(3,3)$
and $\mathcal{M}=\begin{pmatrix}2&1&0\\0&3&5\end{pmatrix}^T$. The \highlight{\smax} losing vectors are given by $(1,2)$ and $(0,4)$, so that Inequality
system~(\ref{ie_feas3}) reads as follows.
$$
  2w_1\ge q,\quad w_1+3w_2\ge q,\quad 5w_2\ge q,\quad w_1+2w_2\le q-1,\,\,\text{ and }\,\, 4w_2\le q-1.
$$
If we minimize one of the objective functions $f_1(q,w_1,w_2)=q$, $f_2(q,w_1,w_2)=w_1$, or $f_3(q,w_1,w_2)=w_2$ $w_1$ subject to those constraints,
we obtain the optimal solution $q=10$, $w_1=5$, $w_2=2$ in all three cases. It is quite remarkable that those values are integers while we have only
requested that they are real-valued. It will turn out that this is a general phenomenon in our context.

Optimal solutions of linear programs are strongly connected with solutions of linear equation systems, since it is well known that, if a linear program admits an optimal solution, then there is an optimal solution attained at a corner of the set of feasible points. To this end we say that an inequality of a linear
program is \textbf{tight} for a given feasible point if \highlight{equality is attained}. In our example the inequalities $2w_1\ge q$, $5w_2\ge q$, and $w_1+2w_2\le q-1$
are tight for the point $(10,5,2)$, while the inequalities $w_1+3w_2\ge q$ and $4w_2\le q-1$ are not. In our context each corner is the solution of an equation system
of three tight inequalities, as we have three (linearly independent) variables.

In Lemma~\ref{lemma_no_three_of_the_same_type}, Lemma~\ref{lemma_2_plus_1}, and Lemma~\ref{lemma_1_plus_2} we check all possible $3$-element subsets
of the inequalities of (\ref{ie_feas3}). It turns out that whenever the corresponding $3\times 3$-equation system has a unique solution, all variables
attain integer values.

So each optimal vertex of the linear program in Lemma~\ref{lemma_reduced_feas_lp} is determined by three tight inequalities of one of the types
$\widetilde{m}^Tw\ge q$ or $\widetilde{l}^Tw\le q-1$, since $w_1,w_2,q\ge 0$ cannot be attained with equality. In the following three lemmas
we consider the possible cases.

\begin{Lemma}
   \label{lemma_no_three_of_the_same_type}
   For Inequality system~(\ref{ie_feas3}), three tight inequalities of type $\widetilde{m}^Tw \ge q$ or three tight
   inequalities of type $\widetilde{l}^Tw\le q-1$ have to be either linearly dependent or do not determine a solution at all.
\end{Lemma}
\begin{Proof}
  Consider the equation system
  $
    aw_1+bw_2=cw_1\highlight{+dw_2=ew_1}+fw_2=z,
  $
  where $z\in\{q,q-1\}$. Eliminating $z$ leaves
  $
    (a-c)w_1+(b-d)w_2=(c-e)w_1+(d-f)w_2=0,
  $
  which has either the unique solution $w_1=w_2=0$, which is infeasible for the whole inequality system, or an infinite number
  of solutions due to scaling. (In the latter case the equations are linearly dependent.)
\end{Proof}

\begin{Lemma}
   \label{lemma_2_plus_1}
   For Inequality system~(\ref{ie_feas3}), two tight inequalities of type $\widetilde{m}^Tw \ge q$ and one tight inequality
   of type $\widetilde{l}^Tw \le q-1$ lead to  an \highlight{integer} solution $(\widehat{q},\widehat{w_1},\widehat{w_2})$ such that
   $w_1\ge\widehat{w_1}$, $w_2\ge\widehat{w_2}$, and $q\ge\widehat{q}$ for all feasible $(q,w_1,w_2)$ or do not determine
   a solution at all.
\end{Lemma}
\begin{Proof}
   Let $(a,b),(c,d)\in\mwvs$ and $(e,f)\in\mlvs$ be the vectors corresponding to the tight inequalities,
   where we assume $a>c$. From $(a,b)\bowtie(c,d)$ and $a>c$ we conclude $d>b+1$. Solving the corresponding equation system yields
   $\widehat{w_1}= \frac{d-b}{Q}$, $\widehat{w_2}= \frac{a-c}{Q}$, and $\widehat{q}= \frac{ad-bc}{Q}$, where 
   $Q:=fc-fa+ad-bc-ed+eb\in\mathbb{Z}$. The case $Q=0$ corresponds to an equation system which does not have a
   unique solution. Since we know that each feasible solution of (\ref{ie_feas3}) satisfies $w_1,w_2>0$ we can assume $Q>0$ in the
   following.
   
   Let $g:=\gcd(a-c,d-b)\ge 1$. For the weights $\widehat{w}_1,\widehat{w}_2$ we can easily check that coalition type $(a,b)$
   has the same weight as coalition type $(a',b')=\left(a-\frac{a-c}{g},b+\frac{d-b}{g}\right)$. If $(a',b')$ is not a winning
   vector, then $(\widehat{q},\widehat{w}_1,\widehat{w}_2)$ cannot be a feasible solution. Thus $(a',b')$ is a \highlight{\smin} winning vector
   too. If $g>1$ then we have $a>a'>c$. We can check $Q':=fc-fa'+a'd-b'c-ed+eb'=\left(1-\frac{1}{g}\right)\cdot Q>0$. Thus we can assume
   w.l.o.g.\ that $a>c$ is minimal within the set of \highlight{\smin} winning vectors corresponding to tight inequalities, i.e.\ we can assume $g=1$.
   
   Now we apply Lemma~\ref{lemma_xgcd} and choose unique integers $u,v$ fulfilling $u(d-b)-v(a-c)=1$, where $0< u\le a-c$ and
   $0\le v<d-b$. The coalition type $(e',f')=(a-u,b+v)$ has weight $\widehat{q}-\frac{1}{Q}$ and thus
   is losing. Since all losing coalitions have weight at most $q-1$ we conclude $Q=1$. Thus $(e',f')$ is indeed a \highlight{\smax}
   losing vector corresponding to a tight inequality. We can easily check $Q'=f'c-f'a+ad-bc-e'd+e'b=1$ so that we can assume $(e,f)=(e',f')$
   since this characterizes the same solution.
   
   Let us have a closer look at the corresponding inequality system again:
   $$aw_1+bw_2-q\ge 0,\quad cw_1+dw_2-q\ge 0,\,\,\text{ and }\,\,-ew_1-fw_2+q\ge 1.$$
   For the basis $(w_1,w_2,q)$ the inverse matrix is given by
   $$
     M^{-1}=\frac{1}{Q}\cdot\begin{pmatrix}d-f&f-b&d-b\\e-c&a-e&a-c\\ed-cf&af-eb&ad-bc\end{pmatrix}.
   $$
   If we can show that all entries of $M^{-1}$ are non-negative, then we have $w_1\ge\widehat{w_1}$, $w_2\ge\widehat{w_2}$, and $q\ge\widehat{q}$ for all feasible
   $(w_1,w_2,q)$.
   
   From $a>c$ and $(a,b)\bowtie (c,d)$ we conclude $a+b<c+d$, so that we have $a-c\ge 1$ and $d-b\ge 2$. Since $e=a-u$, $f=b+v$ with
   $0<u\le a-c$, $0\le v<d-b$ we have $a-e\ge 1$, $f-b\ge0$, $e-c\ge 0$, and $d-f\ge 1$. Thus, the entries of the first two rows of $M^{-1}$ are non-negative
   integers. For $Q=1$ we have $ad-bc=\widehat{q}\ge 1$. From $f=b+v$, $e=a-u$ we conclude $af-eb=av+bu\ge 0$. The last inequality arises from
   $$
     ed-cf\underset{Q=1}{\underbrace{=}}ad-bc-(af-eb)-1=\underset{\ge 1}{\underbrace{a(d-f)}}+b\underset{\ge 0}{\underbrace{(e-c)}}-1\ge0.
   $$
   
   \vspace*{-6mm}
   
   \quad\quad
\end{Proof}

Let us illustrate how Lemma~\ref{lemma_2_plus_1} works by an example. For this purpose let the weighted game $\chi$ be
uniquely characterized by its cardinality vector $\cardv=(4,8)$ and its matrix of \highlight{\smin} winning vectors
$\mathcal{M}=\begin{pmatrix}4&3&2&1&0\\0&1&4&6&8\end{pmatrix}^T$. An {\ireppt} is given by the weights $w_1=7$, $w_2=3$, and
quota $q=24$. Let us assume the that winning vectors $(3,1)$, $(0,8)$ and the losing vector $(1,5)$
would correspond to tight inequalities. The solution of the corresponding equation system is given by
$w_1=\frac{7}{2}$, $w_2=\frac{3}{2}$, $q=12$. Here the weights $w_1$ and $w_2$ are \highlight{non-integer}. So Lemma~\ref{lemma_2_plus_1}
says that $(q,w_1,w_2)$
cannot be a feasible solution of inequality system~(\ref{ie_feas3}). Thus there must be a constraint which is violated. The construction of
$(e',f')$ in the proof precisely gives such a violation. Since $1\cdot (d-b)-2\cdot(a-c)=1$ the coalition $(2,3)$ is a losing vector
with weight $11.5$. We can easily check that it is indeed a \highlight{\smax} losing vector having a weight strictly larger than $q-1=11$.

Starting from the infeasible vector $(12,3.5,1.5)$ the proof provides us even another candidate for a $3$-element subset of tight inequalities.
If we replace the losing vector $(1,5)$ by $(e',f')=(2,3)$, then we obtain the solution $w_1=7$,
$w_2=3$, $q=24$, which now consists of integers. Here we have
$$
 M=\begin{pmatrix}3&1&-1\\0&8&-1\\-2&-3&1\end{pmatrix}\quad\text{and}\quad M^{-1}=\begin{pmatrix}5&2&7\\2&1&3\\16&7&24\end{pmatrix}.
$$
Since the inverse matrix $M^{-1}$ consists of non-negative entries, as generally shown in the proof, we have
$w_1'\ge 7$, $w_2'\ge 3$, and $q'\ge 24$ for every averaged integer representation $(q',w_1',w_2')$. To be more precise: if we combine the
inequalities $3w_1+w_2-q\ge 0$, $8w_2-q\ge 0$, and $-2w_1-3w_2+q\ge 1$ with non-negative multipliers given by the first row on $M^{-1}$, we
conclude $w_1\ge 7$. For the second and third row we similarly obtain $w_2\ge 3$ and $q\ge 24$, respectively. Thus we have found a {\msireppt}.

\begin{Lemma}
  \label{lemma_1_plus_2}
  For Inequality system~(\ref{ie_feas3}), one tight inequality of type $\widetilde{m}^Tw \ge q$ and two tight inequalities
  of type $\widetilde{l}^Tw \le q-1$ lead to an \highlight{integer} solution $(\widehat{w_1},\widehat{w_2},\widehat{q})$ such that
  $w_1\ge\widehat{w_1}$, $w_2\ge\widehat{w_2}$, and $q\ge\widehat{q}$ for all feasible $(w_1,w_2,q)$ or do not determine a solution at all.
\end{Lemma}
\begin{Proof}
    Let $(a,b)\in\mwvs$ and $(c,d),(e,f)\in\mlvs$ be the vectors corresponding to the tight inequalities, where
    we assume $e>c$. Solving the corresponding equation system yields $\widehat{w_1}= \frac{d-f}{Q}$, $\widehat{w_2}= \frac{e-c}{Q}$,
    and $\widehat{q}= \frac{ad-fa+eb-bc}{Q}$, where $Q:=fc-fa+ad-bc-ed+eb\in\mathbb{Z}$. The case $Q=0$ corresponds
    to an equation system which does not have a unique solution. Since we know that each feasible solution of (\ref{ie_feas3}) satisfies
    $w_1,w_2>0$ we can assume $Q>0$ in the following.
   
   Let $g:=\gcd(e-c,d-f)\ge 1$. The vector $(a',b')=\left(e-\frac{e-c}{g},c+\frac{d-f}{g}\right)$ has the same weight as $(a,b)$. So
   similarly to the proof of Lemma~\ref{lemma_2_plus_1} we conclude that $(a',b')$ is a \highlight{\smin} winning vector, which corresponds
   to a tight inequality. We again check that replacing $(a,b)$ by $(a',b')$ is compatible with $Q'>0$ so that we can finally assume $g=1$ w.l.o.g.
   
   Now we apply Lemma~\ref{lemma_xgcd} and choose unique integers $u,v$ fulfilling $u(d-f)-v(e-c)=1$, where $0< u\le e-c$
   and $0\le v<d-f$. The coalition type $(a',b')=(c+u,d-v)$ has weight $q-1+\frac{1}{Q}$. Since losing
   vectors have a weight of at most $q-1$ the vector is winning and we have $Q=1$. Using a similar argument as in the proof of
   Lemma~\ref{lemma_2_plus_1} we conclude that $(a',b')$ is indeed a \highlight{\smin} winning vector corresponding to a tight inequality. We can
   easily check $Q'=fc-fa'+a'd-b'c-ed+eb'=1$ so that we can assume $(a,b)=(a',b')$ since this characterizes the same solution.
   
   Let us have a closer look at the corresponding inequality system again:
   $$
     aw_1+bw_2-q  \ge 0,\quad -cw_1-dw_2+q \ge 1,\,\,\text{ and }\,\,-ew_1-fw_2+q \ge 1.
   $$
   For the basis $(w_1,w_2,q)$ the inverse matrix is given by
   $$
     M^{-1}=\frac{1}{Q}\cdot\begin{pmatrix}d-f&b-f&d-b\\e-c&e-a&a-c\\ed-cf&eb-af&ad-bc\end{pmatrix}.
   $$
   If we can show that all entries of $M^{-1}$ are non-negative, then we have $w_1\ge\widehat{w_1}$, $w_2\ge\widehat{w_2}$, and $q\ge\widehat{q}$ for all feasible
   $(w_1,w_2,q)$.
   
   From $e>c$ and $(c,d)\bowtie (e,f)$ we conclude $e+f<c+d$, so that we have $e-c\ge 1$ and $d-f\ge 2$. Since $a=c+u$, $b=d-v$ with
   $0<u\le e-c$, $0\le v<d-f$ we have $a-c\ge 1$, $d-b\ge0$, $e-a\ge 0$, and $b-f\ge 1$. Thus, the entries of the first two rows of $M^{-1}$ are non-negative
   integers. From $e>c$ we conclude $ed-cf\ge c(d-f)\ge 0$ and from $a=c+u$, $b=d-v$ we conclude $ad-bc=ud+vc\ge 1$. The last inequality arises from
   $$
     eb-af\underset{Q=1}{\underbrace{=}}1+(ed-cf)-(ad-bc)=1+d\underset{\ge 0}{\underbrace{(e-a)}}+c\underset{\ge 1}{\underbrace{(b-f)}}\ge 1.
   $$
   
   \vspace*{-6mm}
   
   \quad\quad
\end{Proof}

\begin{Theorem}
  \label{thm_unique_optimum}
  Let $\chi$ be a weighted game without null voters and with $t=2$, $r>1$. Minimizing the target function $c_1w_1+c_2w_2+c_3q$, where $c_1,c_2,c_3\ge0$
  and $c_1+c_2+c_3>0$, subject to the constraints in (\ref{ie_feas3}) results in a unique optimal integer solution $(q,w_1,w_2)\in\mathbb{N}_{>0}^3$
  satisfying $1\le w_1\le n_2$, $1\le w_2\le n_1$, and $w_2+1\le q\le 2n_1n_2$.
\end{Theorem}
\begin{Proof}
  Let $(q,w_1,q_2)$ be the minimum value of $n_1w_1+n_2w_2$ subject to the constraints in (\ref{ie_feas3}). We already know that the optimum exists.
  This minimum is attained at a corner of the corresponding feasible set and thus arises as the unique solution of a $3\times 3$-equation system,
  consisting of three tight inequalities. Due to Lemma~\ref{lemma_no_three_of_the_same_type} we can apply either Lemma~\ref{lemma_2_plus_1} or
  Lemma~\ref{lemma_1_plus_2}. Thus, each feasible solution $(q',w_1',w_2')$ of inequality system~(\ref{ie_feas3}) has to satisfy $q'\ge q$, $w_1'\ge w_1$,
  and $w_2'\ge w_2$. So we have $c_1w_1'+c_2w_2'+c_3q'\ge c_1w_1+c_2w_2+c_3q$, where equality is attained if and only if $(q',w_1',w_2')=(q,w_1,w_2)$.
  The formulas for $w_1$, $w_2$ and $q$ in Lemma~\ref{lemma_2_plus_1} and Lemma~\ref{lemma_1_plus_2} give the upper bounds $w_1\le n_2$, $w_2\le n_1$,
  and $q\le 2n_1n_2$. Since $\chi$ does not contain null voters we also have $w_1,w_2\ge 1$. If $q\le w_2$, then every single voter would form a
  winning coalition, so that we only have one equivalence class, which contradicts $t=2$.
\end{Proof}

To prove Theorem~\ref{main_thm}, we show that the unique optimal integer solution $(q,w_1,w_2)$ from
Theorem~\ref{thm_unique_optimum} is indeed a {\mirep}. To this end we state that for two feasible solutions $(q,w)$ and $(q',w')$ of Inequality
system~(\ref{ie_feas1}) the vector $\lambda\cdot(q,w)+(1-\lambda)\cdot(q',w')$ is also a feasible solution for all $\lambda\in[0,1]$.

\begin{Lemma}
  Given a weighted game $\chi$ without null voters and with $t=2$, $r>1$, let $(\widehat{q},\widehat{w}_1,\widehat{w}_2)\in\mathbb{N}_{>0}^3$ be
  a feasible solution of (\ref{ie_feas3}), which minimizes the sum of weights $n_1\widehat{w}_1+n_2\widehat{w}_2$. For each {\irep}
  $(q,a_1,\dots,a_{n_1},b_1,\dots b_{n_2})$ of $\chi$ we have $a_i\ge\widehat{w}_1$ for all $1\le i\le n_1$ and $b_j\ge\widehat{w}_2$ for all
  $1\le j\le n_2$.
\end{Lemma}
\begin{Proof}
  It suffices to conclude a contradiction both from $a_{n_1}\le\widehat{w}_1-1$ and $b_{n_2}\le\widehat{w}_2\highlight{-1}$. To shorten the presentation
  we deal with the first case only. Since $\frac{1}{n_1}\cdot \sum_{i=1}^{n_1}a_i\ge\widehat{w}_1$ we can assume $n_1\ge 2$ and 
  \highlight{since every convex combination of a feasible weighting is feasible we can assume} $a_{n_1}=\widehat{w}_1-1$
  w.l.o.g. Next we set $a:=\frac{\sum_{i=1}^{n_1-1}a_i}{n_1-1}$ and $b:=\frac{\sum_{i=1}^{n_2}b_i}{n_2}$. With this the vector
  $$
    (q,\underset{n_1-1}{\underbrace{a,\dots,a}},\widehat{w}_1-1,\underset{n_2}{\underbrace{b,\dots,b}})
  $$
  is also a feasible solution of (\ref{ie_feas1}), where we have $a\ge\highlight{\widehat{w}_1}+\frac{1}{n_1-1}>\highlight{\widehat{w}_1}$.
  
  \highlight{Next we want to utilize the concept of tight inequalities to use a formula between the parameters of the tight inequalities and
  $\widehat{w}_1$. Due to Lemma~\ref{lemma_no_three_of_the_same_type} we have to distinguish the cases of Lemma~\ref{lemma_2_plus_1} and
  Lemma~\ref{lemma_1_plus_2} only.}
  
  If there are two tight inequalities of type $\widetilde{m}^Tw\ge q$ for $(\widehat{q},\highlight{\widehat{w}_1},\highlight{\widehat{w}_2})$,
  see Lemma~\ref{lemma_2_plus_1},
  then let $(c_1,d_1)$, $(c_2,d_2)$ be the two corresponding winning vectors satisfying $c_1>c_2$ and $d_1<d_2$. \highlight{Due to
  Lemma~\ref{lemma_2_plus_1} we have $\widehat{w}_1=d_2-d_1$.} \highlight{Next w}e choose two non-negative integers
  $u\le c_1-c_2$ and $v\le d_2-d_1$ such that $u\cdot(d_2-d_1)-v(c_1-c_2)=1$. We remark $u\ge 1$. With this $(c_3,d_3):=(c_1-u,d_1+v)$ is a losing
  vector corresponding to a tight inequality. Since $(c_1,d_1)$, $(c_2,d_2)$ are winning, $(c_3,d_3)$ is losing,  $c_1\ge 1$, $c_2<n_1$, and $c_3<n_1$
  we have
  $$
    (c_1-1)\cdot a+1\cdot (\highlight{\widehat{w}_1}-1)+d_1\cdot b-q\ge0,\quad
    c_2\cdot a+d_2\cdot b-q\ge0,\,\,\text{ and }\,\,
    -c_3\cdot a-d_3\cdot b+q\ge 1.
  $$
  Summing up $d_2-\highlight{d_3}$ times the first, $d_3-d_1$ times the second, and $d_2-d_1$ times the third inequality yields
  $$
    (\underset{\le 0}{\underbrace{1-d_2+d_3}})\underset{>\widehat{w}_1}{\underbrace{a}}+(\underset{>0}{\underbrace{d_2-d_3}})\widehat{w}_1-
    (\underset{\ge 0}{\underbrace{d_2-d_3}})\ge d_2-d_1,
  $$
  from which we conclude the contradiction
  $$
    \widehat{w}_1>d_2-d_1=\widehat{w}_1.
  $$
  
  If there are two tight inequalities of type $\widetilde{m}^Tw\le q-1$ for $(\widehat{q},\highlight{\widehat{w}_1},\highlight{\widehat{w}_2})$,
  see Lemma~\ref{lemma_1_plus_2}, then let $(c_1,d_1)$, $(c_2,d_2)$ be the two corresponding losing vectors satisfying $c_1<c_2$ and $d_1>d_2$. 
  \highlight{Due to Lemma~\ref{lemma_1_plus_2} we have $\widehat{w}_1=d_2-d_1$.} \highlight{Next w}e choose two non-negative integers
  $u\le c_2-c_1$ and $v\le d_1-d_2$ such that $u\cdot(d_1-d_2)-v(c_2-c_1)=1$. We remark $u\ge 1$. With this $(c_3,d_3):=(c_1+u,d_1-v)$ is a winning
  vector corresponding to a tight inequality. Thus we have 
  $$
    (c_3-1)\cdot a+1\cdot (\widehat{w}_1-1)+d_3\cdot b-q \ge 0,\quad
    -c_1\cdot a-d_1\cdot b+q\ge 1,\,\,\text{ and }\,\,
    -(c_2-1)\cdot a-1\cdot (\widehat{w}_1-1)-d_2\cdot b+q\ge 1.
  $$
  Summing up $d_1-\highlight{d_2}$ times the first, $d_3-d_2$ times the second, and $d_1-d_3$ times the third inequality
  yields the contradiction $\widehat{w}_1>\widehat{w}_1$.
  
  Thus the assumption $a_{n_1}\le \highlight{\widehat{w}_1}-1$ cannot be true and we have $a_i\ge \highlight{\widehat{w}_1}$ for all $1\le i\le n_1$. Similar arguments
  can be outlined for $b_j\ge\widetilde{w}_2$ for all $1\le j\le n_2$.
\end{Proof}
 
\begin{Remark}
  Due to the above lemmas we can algorithmically determine a {\mirep} in
  $O\!\left(\left|\mwvs\right|^3\log(n)\log\log(n)+\left|\mwvs\right|^2\log^2(n)\log\log(n)
  \right)$ time. The case $\left|\mwvs\right|=r=1$ can be dealt directly using Lemma~\ref{thm_min_rep_2_1}. For $r\ge 2$ we consider
  all pairs of \highlight{\smin} winning vectors and all pairs of \highlight{\smax} losing vectors. Here we have
  $\left|\mlvs\right|\le\left|\mwvs\right|+1$ and $\left|\mwvs\right|\le\min\left(n_1+1,
  \left\lfloor\frac{n_2+2}{2}\right\rfloor\right) \le\left\lfloor\frac{n+3}{3}\right\rfloor$ due to Inequality system~(\ref{compact_ilp_2_ge_2}). 
  For each, in $\mwvs\times\mwvs$ or $\mlvs\times \mlvs$ we calculate the parameters $u$ and $v$ via the Euclidean
  algorithm to determine the third tight vector, see Lemma~\ref{lemma_2_plus_1} and Lemma~\ref{lemma_1_plus_2}, respectively. So we have to
  consider at most $\left|\mwvs\right|^2+\left|\mlvs\right|^2$ cases.
  In each case the Euclidean algorithm performs at most $\log(n)$ steps where numbers between $-n$ and $n$ are
  added and divided. After solving the $3\times 3$-equation system, which can be done in time $O(\log(n)\log\log(n))$, we only
  have to check if the solution is feasible.
  Checking the feasibility means determining the minimal weight of a winning vector and the maximal weight of a losing vector,
  which can be done using $O\!\left(\left|\mwvs\right|\right)$ multiplications and additions.
  
  Since the minimal possible values of $w_1$, $w_2$, and $q$ can be bounded via  $w_1\le \max(n_1+1,n_2)$, $w_2\le \max(n_1,n_2-1)$, and
  $q\le (n_1+n_2)\cdot\max(n_1+1,n_2)$ we may also determine a {\mirep} by trying out all possibilities, which results in a
  pseudo-polynomial algorithm.
\end{Remark}

Due to the famous LLL-algorithm \cite{0488.12001,0524.90067} integer linear programs with a fixed number of dimensions, i.e.\ the number of variables,
and a fixed number of constraints can be solved in polynomial time. For a two variables integer program defined by $m$~constraints involving coefficients 
with at most $s$~bits there is a $O(m+\log m\log s)M(s)$ algorithm \cite{1079.90581}, where $M(s)$ is the time needed for $s$-bit integer multiplication
(we assume $M(s)=s\log s\log \log s$). For $t=2$ types of voters we have $\left|\mlvs\right|,\left|\mwvs\right|
\le\left\lfloor\frac{n+6}{3}\right\rfloor$, so that $m=\left|\mwvs\right|\cdot\left|\mlvs\right|+n\in O(n^2)$, and $s\in O(\log n)$ using
the ILP formulation without the quota~$q$. For a general but fixed number of variables Clarkson's sampling algorithm needs an expected number of
$O(m+s\log m)$ arithmetic operations \cite{pre05677125}. Using the ILP formulation with an extra variable for the quota~$q$ we have 
$m=\left|\mwvs\right|+\left|\mlvs\right|+n\in O\!\left(n^{t-1}\right)$ and $s\in O(\log n)$ for $t$~types of voters. We would like to remark that the number of minimal winning vectors can be exponential in $n$ whenever the number $t$ of types of voters is not restricted; see e.g.\ \cite{0841.90134}.

\section{Enumerations and bounds for the number of weighted games}
\label{sec_5}

\noindent
Besides studying properties of complete simple games and weighted games one can also enumerate these special classes of cooperative games
for small numbers of players~$n$. In some cases enumeration results provide a deeper understanding. So far the number of complete simple games of
weighted games is only known up to $n=9$; see e.g.\ \cite{FrMo09,dedekind}. Additionally restricting the parameters $t$ (the number of types of voters)
and/or $r$ (the number of \highlight{\smin} winning vectors) opens the possibility to determine enumeration formulas in some cases. A widely known
result in this context is $csg(n,1)=wvg(n,1)=n$, where $csg(n,t)$ denotes the number of complete simple games with $n$ voters partitioned into $t$
equivalence classes. Similarly $wvg(n,t)$ denotes the number of weighted games with $n$ voters occurring in $t$ different types. In \cite{arxix_freixas}
the authors have determined the formula $cs(n,2)=Fib(n+6)-(n^2+4n+8)$, where $Fib(n)$ denotes the $n$-th Fibonacci number; see
also \cite{dedekind} for an alternative proof. So we know that $cs(n,t)$ is at least exponential in $n$ for $t\ge 2$. In this section we want to show
that the situation changes for weighted games by proving a polynomial upper bound on $wm(n,t)$ in Theorem~\ref{thm_count_1} and Theorem~\ref{thm_count_2}.
It remains to come up with an exact formula for $wm(n,2)$.

If we refine our counts to the numbers $csg(n,t,r)$ and $wvg(n,t,r)$ by additionally considering the number $r$ of \highlight{\smin} winning vectors, more
results can be obtained. In \cite{dedekind} an algorithm is given to principally determine an exact formula for $csg(n,t,r)$ whenever $t$ and $r$
are fixed. So far it is not known whether this can also be done for the number $wvg(n,t,r)$ of weighted games with $t$ types of voters and $r$
\highlight{\smin} winning vectors. For $r=1$ it is not too difficult to come up with such enumeration formulas as we will demonstrate for $t=2$. Having an
exact characterization of the weighted games with $t=2$ and $r=1$ at hand, see the proof of Theorem \ref{thm_min_rep_2_1}, we can easily determine
a formula for their number:

\begin{Corollary}
  \label{cor_wm_2_1}
  The number $wm(n,2,1)$ of weighted games with $t=2$ and $r=1$ is given by $n-1$ for $n\le 2$ and $2(n-2)^2+2$ for $n\ge 3$.
\end{Corollary}

If we skip the parameter $r$ then we can only state an upper bound:

\begin{Theorem}
   \label{thm_count_1}
   $wm(n,2)\le \frac{n^5}{15}+4n^4$.
\end{Theorem}
\begin{Proof}
  Due to the bounds in the minimum integer representation for $r\ge 2$ in Theorem~\ref{main_thm} and Corollary~\ref{cor_wm_2_1} the number $wm(n,2)$ of
  weighted games with $n$ voters and two types of voters is upper bounded by
  $$2(n-2)^2+2+\sum_{n_1=1}^{n-1}\sum_{w_1=1}^{n-n_1}\sum_{w_2=0}^{n_1}\sum_{q=1}^{2n_1(n-n_1)}1=2(n-2)^2+2+2\sum_{n_1=1}^{n-1}(n-n_1)^2(n_1+1)n_1
  \le \frac{n^5}{15}+4n^4.$$
\end{Proof}

For an arbitrary number $t$ of types of voters we can determine the following polynomial upper bound:

\begin{Theorem}
  \label{thm_count_2}
  $$
    wm(n,t)<(tn)^{t^3+2t^2}.
  $$
\end{Theorem}
\begin{Proof}
  Let us denote the weight vector by $w$, the \highlight{\smin} winning vectors by $\widetilde{m}_i$, and the \highlight{\smax}
  losing vectors by $\widetilde{l}_j$. A complete simple game described by the $\widetilde{m}_i$ or the 
  $\widetilde{l}_j$ is weighted if and only if the system of inequalities
  \begin{equation}
    \label{ie_cross_comparision}
    \left(\widetilde{m}_i-\widetilde{l}_j\right)w^T>0   
  \end{equation}
  has a non-negative solution $w$ (for all $i$, $j$). 
  
  Since $\lambda w$ is also a solution for all $\lambda>0$ whenever $w$ is a solution, we consider the equivalent system
  \begin{equation}
    \left(\widetilde{m}_i-\widetilde{l}_j\right)w^T\ge 1.
  \end{equation}
  Such a system of linear inequalities corresponds to a polytope whose vertices correspond to $n$-element subsets of the constraints which are
  attained with equality.  Using the fact that the coefficients of this system of linear inequalities are integers between $-(n-1)$ and $n-1$ we
  can apply Cramers rule to conclude that vertices of this polytope can be written as
  $v_i=\begin{pmatrix}w_1&\dots&w_t\end{pmatrix}=\begin{pmatrix}\frac{a_{2,i}}{b_{2,i}}&\dots&\frac{a_{t,i}}{b_{t,i}}\end{pmatrix}$, where
  $0\le a_{j,i}\le (t-1)!(n-1)^t$ and $1\le b_{j,i}\le (t-1)!(n-1)^t$. Here the common denominator $g$ is bounded from
  above by $\Big((t-1)!(n-1)^t\Big)^{t}$.
  
  Thus multiplying vertex $v_i$ with $g$ yields integer weights $\widetilde{w}_i$ between $0$ and $\Big((t-1)!(n-1)^t\Big)^{t+1}$.
  There are at most $(tn)^{t^3+t^2}$ possible tuples of integer  weights to be considered. The quota can be chosen as the minimum weight
  of a winning coalition. Since there are less than $n^t$ possibilities for the numbers $n_i$ of voters in the $t$ equivalence classes,
  the proposed upper bound on $wm(n,t)$ follows.
\end{Proof}

\section{Concluding remarks}
\label{sec_conclusion}

\noindent
The main result of this paper is that weighted games with two types of players admit a {\mirep}. For three types
of players this need not to be the case. We have shown that by providing examples of games without a {\mirep}. 

We found examples of weighted games with four~types of voters without a {\mireppt}. It is still an open problem
to clarify whether all weighted games with three~types of voters admit a {\mireppt}. To adres this lacuna
we have tried to generalize our technique from Subsection~\ref{subsec_4_2}. One may consider the linear program minimizing the sum of the
weights and have a closer look
at the corners of the corresponding polytope, which are characterized by four equations corresponding to four \textit{tight types of
coalitions} (\highlight{\smax} losing or \highlight{\smin} winning vectors). 

As demonstrated in Subsection~\ref{subsec_4_2} for three tight types of coalitions, the resulting weights and the quota could be
fractional. But using the extended Euclidean algorithm we were able to construct another type of a coalition which contradicts the
tightness of the starting three vectors in these cases. For four tight types of coalitions (and the variables $q$,
$w_1$, $w_2$, and $w_3$) we may go along the same lines and use the extended Euclidean algorithm for three integers in order to deduce some
restrictions on quadruples of tight types of coalitions. This indeed works, but there still remain cases where the optimal
LP solution is fractional. By generating random weighted games with three types of voters we have discovered several such examples, some
of them are given below. For each example we state the sizes of the equivalence classes $\cardv=(n_1,n_2,n_3)$, the 
\highlight{non-integer} minimum sum representation preserving types $\tau_r=(q,w_1,w_2,w_3)$, and the {\msireppt}
$\tau_i=(q,w_1,w_2,w_3)$:
\begin{itemize}
 \item[(1)] $\cardv=(9,62,71)$, $\tau_r=(154.\overline{3},38.\overline{3},22.\overline{6},6.\overline{6})$,
            $\tau_i=(185,46,27,8)$
 \item[(2)] $\cardv=(19,52,65)$, $\tau_r=(3984.2,200,110,76.6)$,
            $\tau_i=(5617,282,155,108)$            
 \item[(3)] $\cardv=(30,93,30)$, $\tau_r=(122.\overline{3},22.\overline{3},16,9.\overline{3})$,
            $\tau_i=(131,24,17,10)$
 \item[(4)] $\cardv=(8,99,10)$, $\tau_r=(51,17,10.5,4.5)$,
            $\tau_i=(57,19,12,5)$
 \item[(5)] $\cardv=(3,71,37)$, $\tau_r=(347.5,100,31.5,15)$,
            $\tau_i=(441,127,40,19)$           
\end{itemize}
Originally we have obtained the values of $\tau_i$ by minimizing $n_1w_1+n_2w_2+n_3w_3$ but it turned out that in all
of these (and the other found) cases we have a {\mireppt}, so that minimizing $w_1$, $w_2$, $w_3$, or $q$ would yield the same result.
We would like to remark that we have also found some example where only one value is \highlight{non-integer}. Although in our experiments the only occurring
denominators were $2$, $3$, and $5$, we do not think that the denominators are bounded by a constant. So far we have a very poor
probabilistic model which generates those examples with a very low probability. Nevertheless we have a strong feeling that each weighted game with
three types of voters admits a {\mireppt}. As a small justification we would like to remark that we have tried some specific parametric constructions
which provably do not contain counter examples. 

We leave the challenging question of whether each weighted games with three types of voters admits an {\mireppt}
open for the interested reader and hope that our specific examples might help to get some
useful insights. One can get a first glimpse of the difficulty of this problem by comparing the values of $\tau_r$ and $\tau_i$
in our examples.

\medskip

Weighted games with an arbitrary number of {\msirep s} have been generated in Subsection \ref{subsec_3_2}. Moreover, some bounds have been
obtained for the number of non-isomorphic weighted games depending on the number of voters and on the number of types of voters, and the existence of
a weighted game, in {\mirep} for any pair of two coprime integer weights, has been determined.

Other interesting open problems in the context of this paper are the question for a weighted game with a unique {\msirep}, but without a
{\mirep}, and the question for a polynomial time algorithm to determine {\msirep s} for weighted games or a proof that this problem is $NP$-hard.

Another important line of research would be to deepen our understanding of the link between {\mirep s} of weighted games and one-point 
solution concepts, like the nucleolus, least core, etc.; see e.g.\ \cite{0841.90134,peleg}.

Of course the techniques presented in this paper may be applied to study similar questions for roughly weighted games.

\section*{Acknowledgements}

\noindent
The authors thank Stefan Napel for carefully reading a preliminary version of this article. The research of the first author
was supported by Grant SGR 2009-01029 of The Catalonia Government (\textit{Generalitat de Catalunya}) and Grant MTM 2012-34426 
of the Spanish Economy and Competitiveness Ministry. He also acknowledges the Barcelona Graduate School of Economics and the
Generalitat of Catalunya for their support.


\end{document}